\documentclass[a4paper,oneside,12pt]{article}
\usepackage{amsmath,amsfonts,amscd,amssymb}
\usepackage{longtable,geometry}
\usepackage{graphicx}
\usepackage[english]{babel}
\usepackage[utf8]{inputenc}
\usepackage[T1]{fontenc}
\newtheorem{theorem}{Theorem}[section]

\newtheorem{lemma}[theorem]{Lemma}
\newtheorem{proposition}[theorem]{Proposition}

\numberwithin{equation}{section}

\newcommand{\IP}{\mathbb{P}}
\newcommand{\IR}{\mathbb{R}}
\newcommand{\IE}{\mathbb{E}}
\newcommand{\B}{\mathcal{B}}
\newcommand{\Zd}{\mathbb{Z}^d}

\title{The Limiting Shape for Drifted Internal Diffusion Limited Aggregation is a True Heat Ball}
\author{Cyrille Lucas\footnote{Mod\'elisation al\'eatoire de Paris 10 (MODAL'X) \newline \textit{email:} cyrille.lucas@u-paris10.fr}}
\date{28/01/2012}
\begin{document}

\maketitle


\begin{center}
\small\textsc{Abstract:} We build the iDLA cluster using drifted random walks, and study the limiting shapes they exhibit, with the help of sandpile models. For constant drift, 
 the normalised cluster converges to a canonical shape $S$, which can be termed a true heat ball, in that it gives rise to a mean value property for caloric functions. The existence and boundedness of such a shape answers the natural yet open question of the existence and boundedness of a shape that satisfies a mean value property for caloric functions.
 \normalsize

\bigskip

\textbf{Keywords:} Internal Diffusion Limited Aggregation, Limiting Shape, Parabolic Free boundary Problem, Heat Equation\footnote{\textbf{MSC classes:} 60G50, 35K05, 35R35.}.

\end{center}

\newpage

\tableofcontents

\section{Introduction}

The internal Diffusion Limited Aggregation (iDLA) model was first introduced by Diaconis and Fulton in \cite{DF} and gives a protocol for building a random set recursively. At each step, the first vertex visited outside the cluster by a random walk started at the origin is added to the cluster.

The question of the limiting shape of this model is a well-studied one. Lawler, Bramson and Griffeath \cite{lawler1992internal} were the first to identify the limiting shape of the model, in the case of simple symmetric random walks, as a Euclidean ball. Their result was later sharpened by Lawler \cite{lawler1995subdiffusive}, who gave a polynomial upper bound for the fluctuations of the aggregate around this limiting shape. The question of the fluctuations of the aggregate around the limiting shape recently became of renewed interest with the simultaneous works of Asselah and Gaudill\`ere (\cite{asselah2010logarithmic}, \cite{asselah2010sub}) and Jerison, Levine and Sheffield (\cite{jerison2010internal}, \cite{jerison2010logarithmic}), who provided a logarithmic upper bound for the fluctuations. Jerison, Levine and Sheffield have since provided a partial description of the fluctuations that relates them to the Gaussian Free Field (\cite{jerison2011internal}). 

Still using simple random walks, Peres and Levine introduced in \cite{levine2010scaling} a method for identifying limiting shapes of a larger class of iDLA models, where all particles are not started at the origin, but rather from multiple starting points or even from an initial density. They proved convergence towards limiting shapes that are not Euclidean balls.

All these results are proved in the framework of simple random walks, and other random walks did not appear in the literature until Blach\`ere's article \cite{blachere2002agregation}. In this paper, the iDLA model is studied for centered random walks, and convergence towards the ball of a specified norm is proved under moment conditions on the random walks. At the end of the paper, the case of drifted random walks is studied, and a limiting shape is found in the one-dimensional case. The author then conjectures the existence of a limiting shape for all dimensions. The initial idea that the limiting shape should be a level line of the Green's function, as it is the case in several types of groups (see \cite{blachere2007internal}, \cite{blachere2004internal}) is disproved.

In this paper, we present a limiting shape result for a simple class of drifted random walks. The limiting shape of the normalised cluster is characterised as a true heat ball because it gives rise to a mean-value property for caloric functions. The existence of such a bounded shape is an open problem in PDE theory (see \cite{hakobyangeneralized}), for which our convergence provides a proof inherited from the field of random walks.

For the sake of simplicity, we will consider the following class of drifted random walks (Our result will be extended to a more natural class of drifted random walks in section \ref{extension}). For $p \in (0,1)$, let $ (S^j)_{j \in \mathbb{N}}$ a sequence of independent random walks on $\Zd$, with the following law:
\begin{eqnarray}
\IP\left( S (t+1) - S(t) = \pm e_i \right) &=& \frac{1-p}{2(d-1)} \qquad \text{ for $i=1 \cdots d-1,$ and} \nonumber\\
\IP\left( S (t+1) - S(t) = e_d \right) &=& p. \nonumber\\ 
\end{eqnarray}
We will build our cluster using the sequence of random walks $(S^j)_{j \in \mathbb{N}}$, then normalise it as specified in section \ref{normalisation}. Specifically, we will normalise a point $x=(x_1, \cdots, x_d) \in \Zd$ to a point $(z,t) \in \IR^{d-1}\times\IR$. In order to simplify notations, we call the first $d-1$ coordinates \emph{space coordinates} (denoted by $z$) and the drift coordinate the \emph{time coordinate} (denoted by $t$).

The main result in this paper is the following:
\begin{theorem}
\label{maxitheoreme}
Let $\mathcal{A}_n$ be the normalised drifted internal diffusion limited aggregation cluster. Then there exists a set $D \subset \IR^{d-1}\times\IR_+$ with the following properties:
\begin{enumerate}
\item  Almost surely, $\mathcal{A}_n$ converges towards $D$ with respect to the Hausdorff distance.
\item Let $\phi$ be a $C^{\infty}$ function of time and space such that
$$ \frac{1-p}{2(d-1)} \Delta \phi + p \frac{\partial \phi}{ \partial t} =0. $$
Then the following mean value property holds:
$$ \int_{D} \phi(z,t) dz dt = |D| \phi(0). $$
\item The set $D$ is bounded in time and space.

\end{enumerate}
\end{theorem}

\begin{figure}[ht]
\centering
\includegraphics[scale=0.65]{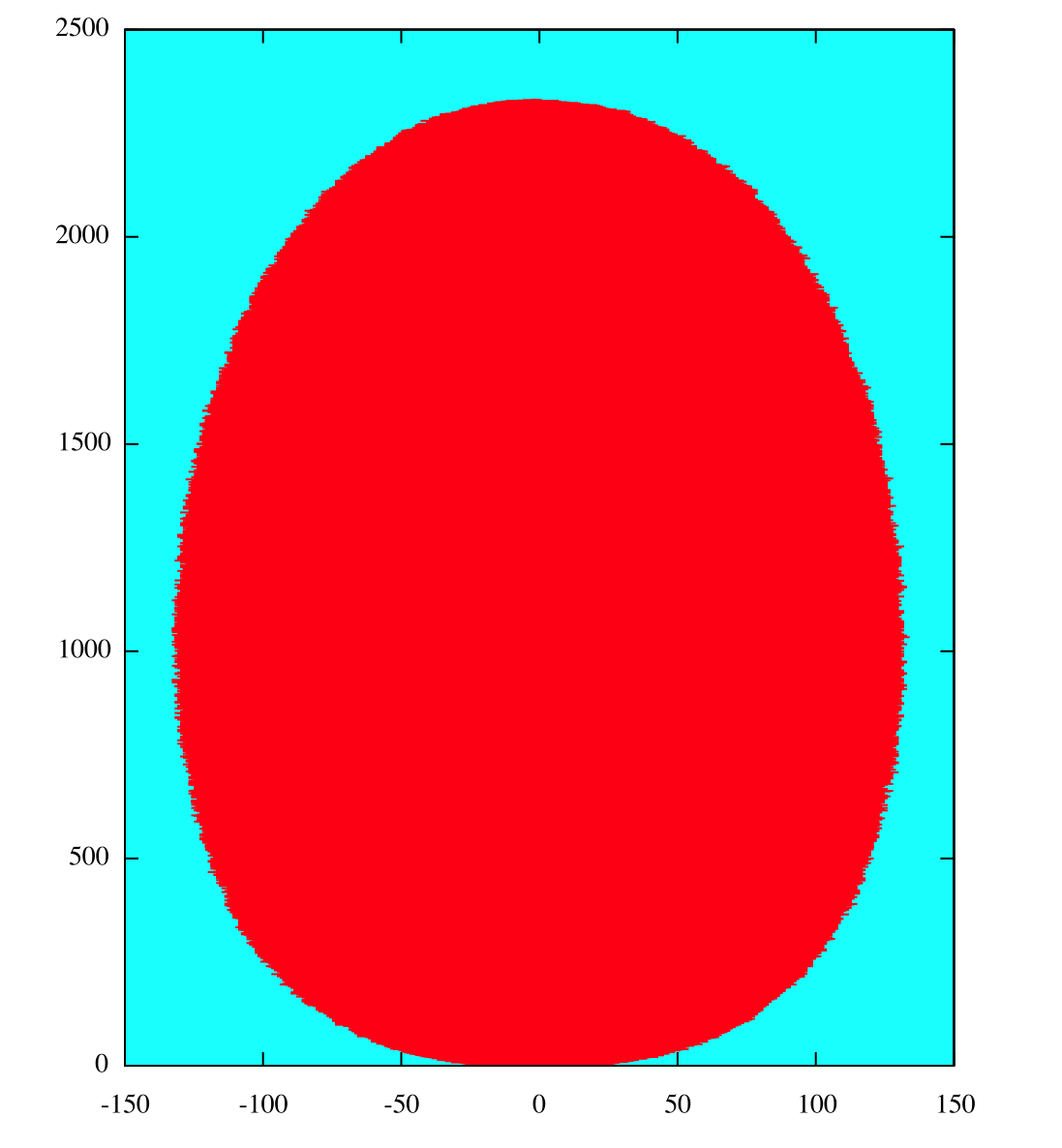}
\caption{\textit{Drifted iDLA aggregate with 500 000 particles and  $p = 0.2$ }}\label{fig:1}
\end{figure}

The first point of the theorem states our convergence result for the normalised aggregate. The model admits a deterministic limiting shape whose properties are detailed in points 2 and 3. The second point states a mean value property for caloric functions, using the shape $D$, which justifies our use of the term true heat ball. The last point provides an answer to the as yet open problem of finding a such a mean value property on a bounded set. In a nutshell, the normalised iDLA cluster converges almost surely to a bounded true heat ball. Figure \ref{fig:1} is a simulation of this shape, obtained with 500 000 particles.

Our proof follows the general idea of Levine and Peres, whose method for finding limiting shapes can be translated in our context. First, we introduce an equivalent of the divisible sandpile model (for a definition and convergence result of the original model, see \cite{levine2009strong}), which we call the unfair divisible sandpile, because one direction is privileged throughout the construction of the cluster.
We study this model in details, and provide a limiting shape result. We prove convergence towards an abstract shape $D$, which for now is defined through a parabolic obstacle problem.

Parabolic obstacle problems are relatively frequent in the literature, and occur in a number of varied fields. First and foremost, they are studied in the context of heat diffusion, and in particular in the case of two-phase transition equations, like the Stefan problem (see \cite{meirmanov1986stefan} and \cite{crank1987free} for a discussion of general parabolic free boundary problems). They also appear in finance, namely in problems related to the pricing of American put options: see for example \cite{kuske1998optimal}, \cite{pham1997optimal}, \cite{jacka1991optimal}. In our case, we will use very strong results in parabolic potential theory (see for example \cite{caffarelli2004regularity}, \cite{athanasopoulos1996caloric}) to ensure that the limiting shape $D$ is smooth enough, and to characterise it as a true heat ball.

The next step of our proof is to relate the drifted iDLA and unfair divisible sandpile, and prove that they almost surely share the same limiting shape $D$. To conclude, we use probabilistic arguments to give bounds on the drifted iDLA cluster, which are in turn used to bound $D$, thus proving the convergence of our normalised cluster towards a bounded true heat ball.


\section{Heuristics}

In this section, we will motivate the introduction of the divisible sandpile model for the study of the iDLA model. We consider the \textit{odometer function} introduced in \cite{levine2010scaling}, that measures the total number of particles emitted from point $x \in \Zd$ throughout the construction of the cluster, counted with repetitions.

A given point will then start with a certain number of particles ($n$ if it is the origin, 0 otherwise, in our situation), receive new particles during the construction of the cluster, and end with a new number of particles (1 if it eventually lies inside the cluster, 0 otherwise).

Of all the particles that passed through a given point $x$, it seems natural that a proportion travelled to each of the reachable neighbors of $x$, and that proportion should in some sense be close to the transition probability from $x$ to that particular neighbor.

Assuming that the same holds for each of the neighbors of $x$, we get a tentative local equation:
$$\sum_{y \sim x} p(y,x) u(y) - u(x) \stackrel{?}{=}  \nu(x) - \sigma(x), \label{heuristic1}$$
where $\sigma(x)$ and $\nu(x)$ are the initial and final amounts of mass at point $x$, respectively.

While this reasoning is flawed because it assumes independence of correlated quantities (among other problems), it motivates the introduction of a new model in which such a local relation holds. The unfair divisible sandpile model will play this role in our case.



\section{Unfair divisible sandpile}
\subsection{Definitions and notations}

We introduce the unfair divisible sandpile model to be the drifted counterpart of the divisible sandpile model defined in \cite{levine2009strong}. Consider a continuous distribution of mass on $\Zd$, with finite total mass and bounded support. A lattice site is full if it has mass at least $1$. Any full site can topple by keeping mass 1 for itself, and distributing the excess mass among its neighbors. While in the classical divisible sandpile model, the mass is split equally among neighbors, in our model, it will be distributed proportionally to the step distribution of the drifted random walk $S$.

At each time step, a full site is toppled. When we let the time go to infinity, if every full site is toppled infinitely often, the mass converges to a limiting distribution in which each site has mass less than $1$. This is the object of lemma \ref{convergenceDS}.

Just like the divisible sandpile model, our unfair version is abelian, in the following sense: while individual toppling do not commute, the limiting distribution of mass does not depend on the order of the topplings, provided that they are done in an appropriate fashion, which is to say that each site that becomes full in the course of the topplings is then toppled infinitely often.

We prove this abelian property in lemma \ref{abelian}.

A crucial tool for the study of the divisible sandpile model, be it classical or unfair, is the odometer function. We define it as the mass emitted from $x$ throughout the construction of the cluster:
$$u(x) = \text{total mass emitted from } x.$$

It is important to note that this quantity does not depend on the sequence of topplings either, which will derive from our proof of lemma \ref{abelian}. Since the mass emitted from a point is always distributed in the same fashion, we can remark that the mass received by $x$ from his neighbors is the following:
$$ \text{mass received in $x$} =p u(x - e_d) + \frac{1-p}{2(d-1)} \sum_{y \sim x, y-x \bot e_d} u(y) $$
Since we have defined $u$ as the emitted mass, the difference between the received mass and $u$ at point $x$ will be equal to the difference between the initial and final amounts of mass at point $x$. Namely,
$$ \nu(x) - \sigma(x) =p u(x - e_d) + \frac{1-p}{2(d-1)} \sum_{y \sim x, y-x \bot e_d} u(y) - u(x),  $$
where $\sigma$ and $\nu$ are the initial and final amounts of mass. We split the quantity $u$ so as to make sense of this equation:
\begin{eqnarray}
\nu(x) - \sigma(x) &=& p u(x - e_d) + \frac{1-p}{2(d-1)} \sum_{y \sim x, y-x \bot e_d} u(y) - u(x) \nonumber\\
 &=& p \left(u(x - e_d) - u(x) \right) + \frac{1-p}{2(d-1)} \sum_{y \sim x, y-x \bot e_d} \left(u(y) - u(x)\right) \nonumber\\
 &=&  - p \left( u(x) - u(x-e_d)\right) + (1-p) \tilde{\Delta} u(x), \label{unimportant1}
\end{eqnarray}
where $\tilde{\Delta}$ is the discrete Laplacian taken over the first $d-1$ coordinates of $\Zd$, which is to say that: $ \tilde{\Delta} u(x) = \frac{1}{2(d-1)} \sum_{y \sim x, y-x \bot e_d} \left(u(y) - u(x)\right)$.

We now define a discrete operator $\mathcal{K}$ that sums up this operation, and will play an important role in our proofs:
$$ \mathcal{K}f(x) =  (1-p) \tilde{\Delta} f(x)  - p \left( f(x) - f(x-e_d)\right). $$

We can now restate equation (\ref{unimportant1}) in the following way:
$$\nu(x) - \sigma(x) = \mathcal{K}u (x) \label{important1}$$

Given the nature of $\mathcal{K}$, and since we will be inclined to consider its continuous counterpart, we will transform notations when we pass to the continuous: the first $d-1$ will correspond to the space coordinates, while the last coordinate will correspond to the time coordinate $t$. The continuous counterpart of $\mathcal{K}$ will be the following operator:
$$ \mathfrak{K} f(x,t) = \frac{1-p}{2(d-1)} \Delta f(x,t) - p \frac{\partial f}{\partial t} (x,t),$$
which is well known as the heat operator. Hence, we will henceforth call $\mathcal{K}$ the \textit{discrete heat operator}, and define a discrete caloric (respectively supercaloric) function as a function $f: \Zd \rightarrow \IR$ such that $\forall z \in \Zd, \mathcal{K}f(z) = 0$ (respectively $\mathcal{K}f(z) \leq 0$).

The sense in which $\mathfrak{K}$ will be the continuous counterpart of $\mathcal{K}$ will be a non-trivial point in our proof, and will give rise to the non-standard normalisation we introduce in section \ref{normalisation}.

\subsection{Convergence and abelian property}
\label{nestintermediaire}

We will now state the convergence result and abelian property for the unfair divisible sandpile model. Let us first define a \textit{toppling scheme}. A \textit{toppling scheme} T for a configuration $\nu_0$ is an infinite sequence of indexes in $\Zd$ in which each site that is initially full or becomes full through the realisation of previous topplings appears infinitely often in the remainder of the sequence. Note that we exclude from consideration the schemes that terminate after a finite number of topplings.
We will denote by $\nu_k^T$ the distribution of mass after the toppling of the first $k$ sites listed in $T$ (in the specified order), and $u_k^T(x)$ the mass emitted from $x$ up to the $k$-th toppling. We omit superscripts when only one toppling scheme is involved.

In addition, a toppling scheme will be called \textit{legal} if it only tries to topple full sites.

\begin{lemma}
\label{convergenceDS}
Consider an initial configuration $\nu_0$ with a finite total amount of mass $M$ and bounded support $\mathcal{S}$, and a legal toppling scheme $T$ for this distribution. Then as $k$ tends to infinity, $u_k$ and $\nu_k$ tend to limits $u$ and $\nu$. The limiting configuration $\nu$ is such that $\forall x \in \Zd, \nu(x) \leq 1$. Moreover, these limits satisfy the following mass relation:
$$ \nu = \nu_0 + \mathcal{K}(u). \label{mass1}$$
\end{lemma}

Proof: Let $B$ be a ball centered at the origin big enough to contain all points within graph distance $M$ of $\mathcal{S}$. We remark that any point that if $\nu_k(x)>0$ at one point in the construction of the cluster, it means either that $x$ was initially in the support of $\nu_0$ or that it has received mass from a neighbor with mass greater that one. Since the same reasoning can be applied to this particular neighbor, it means that all points that have positive mass must be within graph distance $M$ of $\mathcal{S}$. Hence, $\nu_k$ is supported in $B$.

We introduce the weight function:
$$W_k = \sum_{x \in \Zd} \nu_k(x) \left( x_1^2+ \cdots +x_{d-1}^2 + x_d \right).$$
Note that $W$ can be negative, which does not present any problem as the key property is only that it increases through topplings.

When toppling site $x$ at time $k$, the mass at point $x$ has been modified by an amount $\alpha_k(x) =\nu_{k-1}(x)-\nu_{k}(x)$ which has been transferred to the neighbors of $x$ according to the toppling rule. Hence,
\begin{eqnarray}
W_k - W_{k-1} &=& \alpha_k(x) \frac{1-p}{2(d-1)}\sum_{y \sim x, y-x \bot e_d} y_1^2-x_1^2+ \cdots +y_{d-1}^2-x_{d-1}^2 \nonumber\\
& & \qquad \qquad \qquad \qquad \qquad \qquad + p \alpha_k(x) \left(x_d+1 - x_d \right)  \nonumber\\
&=& \alpha_k(x)
\end{eqnarray}

Since $u_k$ is the sum up to $k$ of all the relevant $\alpha_i(x)$, we get the following relation on weights:
$$ W_k = W_0 + \sum_{x\in \Zd} u_k(x).$$
For every $x$, $u_k(x)$ is a bounded increasing sequence, so it converges to a value $u(x)$.

At all finite times $k$, the relation
$$\nu_k(x) = \nu_0(x) + \mathcal{K} u_k(x) \label{mass2}$$
holds, so that the convergence of $\nu_k$ is a consequence of that of $u_k$. Moreover, its limit $\nu$ is such that $\nu = \nu_0 + \mathcal{K}u$.

Now a point $x$ is either never toppled, in which case we have $\nu_k(x) \leq 1$ for all values of $k$, or it is toppled infinitely often, but then each time a toppling occurs at point $x$, $\nu(x) \leq 1$. In both cases, the limit $\nu(x)$ has to satisfy $\nu(x) \leq 1$, which concludes the proof.

We next prove the abelian property of the unfair divisible sandpile model.
\begin{lemma}
\label{abelian}
Consider an initial configuration $\nu_0$, and two legal toppling schemes $S$ and $T$. Then the corresponding final configurations $\nu^{S}$ and $\nu^{T}$ are equal, as are the final odometer functions $u^{S}$ and $u^{T}$.
\end{lemma}
Proof: Let $x_k$ be the point at which the $k$-th toppling occurs in the scheme $S$. We will prove by induction on $k$ that $u^{T} (x_k) \geq u_k^{S}(x_k)$. Note that we are comparing the final version of the odometer $u^T$ to the partially toppled odometer $u^S_k$.

The base case of the induction is $u^T(x_1) \geq u_1(x_1)= \sigma(x_1)-1$, which is always true, because $x_1$ is eventually toppled in scheme $T$.

Suppose that the property holds for all $i < k$. Then for $x \neq x_k$, either $x$ is not toppled before time $k$ in scheme $S$, and $u_k^S(x) = 0$, or it is, in which case we consider the last index $i$ for which toppling occurs at point $x$. Then $u^T(x) \geq u_i^S(x) = u_k^S(x)$. In both cases, $u^T(x) \geq u_k^S(x)$.

Since $S$ is legal, and $T$ is a toppling scheme,
$$ \nu^T(x_k) \leq 1 \leq \nu_k^S (x_k). $$
Hence, equations (\ref{mass1}) and (\ref{mass2}) ensure that:
$$\mathcal{K} u^T (x_k) \leq \mathcal{K} u^S (x_k).$$

We write the mass relations for point $x_k$:
\begin{eqnarray}
u^T(x_k) - \mathcal{K}u^T(x_k) = \nu_0 (x_k) - \nu^T(x_k) \nonumber\\
u^S_k(x_k) - \mathcal{K} u^S_k(x_k) = \nu_0 (x_k) - \nu^S_k(x_k). \nonumber
\end{eqnarray}

Taking the difference yields:
\begin{eqnarray}
u^T(x_k)- u^S_k (x_k) &\geq& {\frac{1-p}{2(d-1)} \sum_{y \sim x_k, y-x_k \bot e_d} \left( u^T(y) - u_k^S(y) \right) \break  }\nonumber\\
 & & \qquad \qquad + {p \left( u^T(x_k-e_d) - u^S_k(x_k-e_d) \right)} \nonumber\\
 &\geq & 0. \nonumber
\end{eqnarray}
The right hand side is indeed positive since it only involves differences for points different from $x_k$.

Hence, we have proved by induction that for all $k$, $u^{T} (x_k) \geq u_k^{S}(x_k)$. It follows that $u^T \geq u^S$. A symmetric argument shows that $u^S \geq u^T$, which concludes the proof.

\subsection{Limiting shape}
\subsubsection{Normalisation}
\label{normalisation}
Since the normalisation we will use is non-standard, we introduce a specific normalisation function $\phi_n$. The fact that this normalisation is not the same in the direction of the drift as in the other directions plays a very important part throughout the paper. It is defined as follows, for functions defined on $\Zd$ with values in $\IR$ :
$$\begin{array}{lrcl}
\forall f: \Zd \rightarrow \IR, & \phi_n (f) & : & \IR^{d-1} \times \IR \rightarrow \IR \\
 & \phi_n(f)(x,t) &=& f\left( (\lfloor n^{\frac{1}{d+1}} (x)_i \rfloor)_{i \in \{1,...,d-1\}} , \lfloor n^{\frac{2}{d+1}} t \rfloor \right).
\end{array}$$
With a slight abuse of notation, we extend the definition of $\phi_n$ to subsets $A$ of $\Zd$:
\begin{eqnarray}
\phi_n (A) &=& \left\{ (x,t) \in  \IR^{d-1} \times \IR, \left( (\lfloor n^{\frac{1}{d+1}} (x)_i \rfloor)_{i \in \{1,...,d-1\}} , \lfloor n^{\frac{2}{d+1}} t \rfloor \right) \in A \right\}. \nonumber
\end{eqnarray}


\subsubsection{The heat equation}
\label{heatequation}
In the following sections, we will be required to handle the heat equation, and in particular to define what a supercaloric function is in the continuous setting. We will follow exactly the definition of Evans, \cite{Evans}. First, we define the fundamental solution $\Phi(x,t)$ of the heat equation $\mathfrak{K}(f)=0$ as follows:
\begin{eqnarray*}
 \Phi(x,t) = \begin{cases} \left(\frac{\beta}{\pi t} \right)^{\frac{d-1}{2}} \mathrm{e}^{-\beta \frac{x^2}{t}} & \text{for $t>0,$}
\\
0 &\text{for $t<0,$}
\\
0 &\text{for $t=0, x \neq 0$.}
\end{cases}
\end{eqnarray*}
where we define
$$ \beta = \frac{p(d-1)}{2(1-p)}.$$
Note that $\Phi$ is singular at the origin, and solves $\mathfrak{K}(\Phi)=0$ away from the origin.

For fixed $x\in \IR^{d-1},t \in \IR, r >0$, we define:
$$E_{(x,t;r)} = \left\{ (y,s) \in \IR^{d-1} \times \IR, s \leq t, \Phi(x-y, t-s) \geq \frac{1}{r^{d-1}} \right\}.$$

Note that the boundary of $E$ is a level set of $\Phi(x-y, t-s)$. Point $(x,t)$ is at the top, center end of $E$.

Let $f$ be a measurable function; we define the operator $\mathfrak{A}_r$ as follows:
$$ \mathfrak{A}_r(f)(x,t) = \frac{ \beta }{4r^{d-1}} \int_{E_{(x,t;r)}} f(y,s) \frac{|x-y|^2}{(t-s)^2} dy ds.$$

This operator is a sort of mean value operator; indeed, if $f$ is a smooth solution of the heat equation, $\mathfrak{A}_r$ can be used to compute the value at a given point:
\begin{eqnarray}
u(x,t) = \mathfrak{A}_r (u)(x,t). \label{meanVP}
\end{eqnarray}

Note that the right hand side involves only $u(y,s)$ for times $s \leq t$, which is reasonable, because the value at a given time should not depend upon future values of $u$.

We now use this operator to define supercaloric functions, as follows. Let $f$ be a measurable function. We say that $f$ is supercaloric if it satisfies the following property:
$$ \forall (x,t) \in \IR^{d-1} \times \IR, \forall r>0, \qquad f(x,t)  \geq \mathfrak{A}_r (f)(x,t). $$

For a function $\phi$ that has $C^{2,1}$ regularity, this is equivalent to the intuitive result:
\begin{equation*}
 \mathfrak{K}\phi \leq 0.
\end{equation*}

Note that equation (\ref{meanVP}) is only a weak substitute for the mean-value property of harmonic functions; instead of simply integrating the function on a given shape, a kernel is used. The problem of the existence of a generalized mean value property for caloric functions, which is a consequence of Theorem \ref{maxitheoreme}, is discussed in \cite{hakobyangeneralized}.

\subsubsection{Main result}

The following theorem provides a result for the limiting shape of the unfair divisible sandpile aggregate under the proper normalisation.

We run the unfair divisible sandpile on $\Zd$ with the following initial mass configuration:
$$ \sigma_n(x) = n\delta_0(x), $$
and we intuitively define the unfair divisible sandpile aggregate as the set of points in $\Zd$ that have positive final mass:
$$ D_n = \{ x \in \Zd, \nu_n(x) =1 \}, $$
where $\nu_n$ is the final mass configuration corresponding to the inital mass distribution $\sigma_n$. Conversely, $u_n$ is the odometer function corresponding to the inital mass distribution $\sigma_n$. Note that this is a new notation and replaces the notation introduced in section \ref{nestintermediaire}. The functions $u_n$ and $\nu_n$ are the (unique) values obtained as the limits of any legal toppling scheme.

\begin{figure}
\includegraphics[scale=0.18]{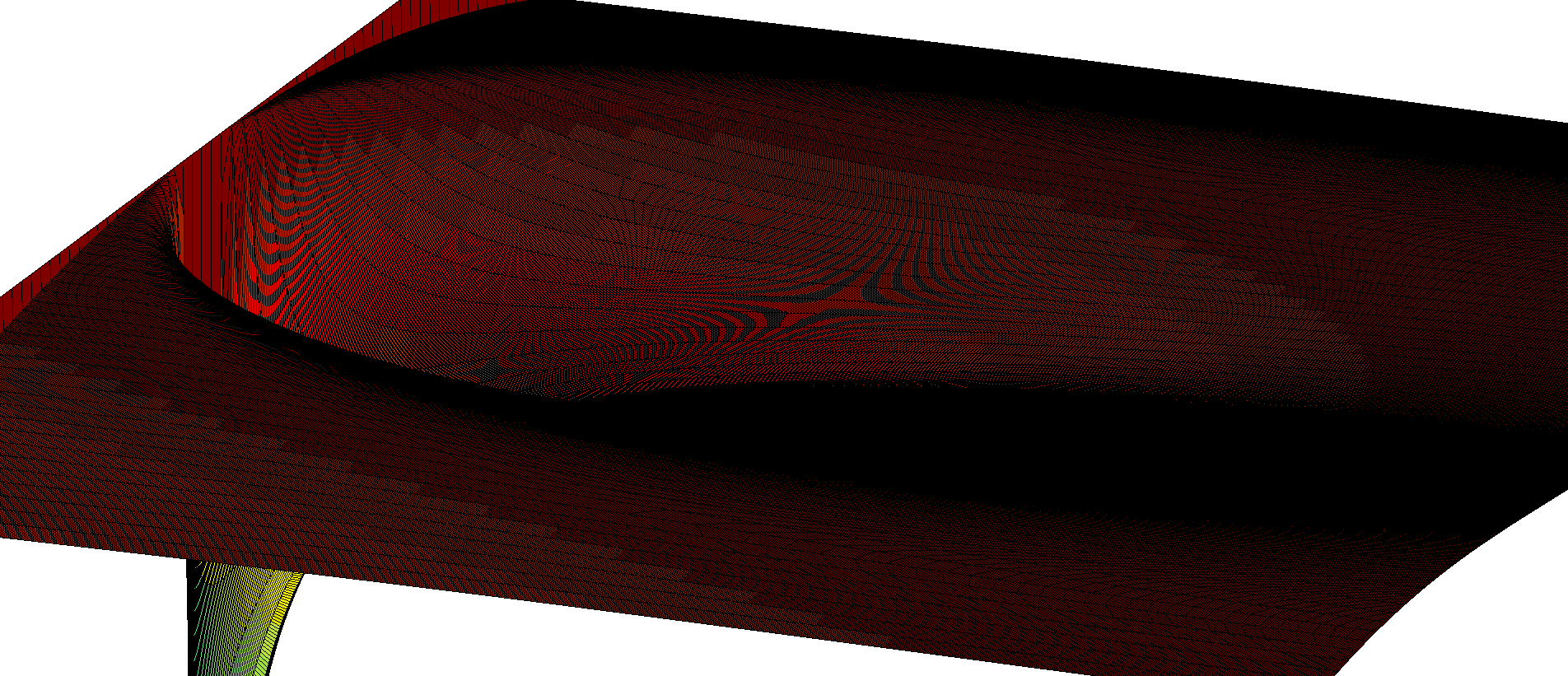}
\caption{The obstacle function $\gamma(x,t)$.}\label{fig:2}
\end{figure}

Recall that $ \beta  = \frac{p(d-1)}{2(1-p)},$ and define the following quantities:
\begin{eqnarray}
\gamma (x,t) & = &   t - ||x||^2  - \frac{1}{p} \left( \frac{\beta}{ \pi t}\right)^{\frac{d-1}{2}} \exp\left( -  \beta \frac{||x||^2}{t} \right) ,\nonumber\\
s      (x,t)& = & \inf\{h(x,t)|  h(x,t)\text{ is supercaloric and }h\geq \gamma \}, \nonumber\\
u      (x,t)& = & s(x,t)- \gamma(x,t) \nonumber\\
D      & = & \{ (x,t) \in \IR^{d-1} \times \IR^*_+, u(x,t)>0 \} = \{ s- \gamma >0 \}. \nonumber
\end{eqnarray}
The function $\gamma$ is plotted in Figure \ref{fig:2} when $d=2$.

Then, the following theorem states the convergence of the normalised unfair divisible sandpile aggregate on compacts.
\begin{theorem}
For every compact $K$ of  $\IR^{d-1}  \times \IR^*_+ $, the intersection of $K$ and the normalised unfair divisible sandpile aggregate $\phi_n(D_n) \cap K$ converges to $D\cap K$ with respect to the Hausdorff distance. 
\label{sandpile}
\end{theorem}


We start by proving that the odometer function for the unfair divisible sandpile with starting configuration $n$ at the origin is given by the difference between an obstacle function $\gamma_n$ (for which we fix $\mathcal{K}\gamma_n$) and its least supercaloric majorant.
\begin{lemma}
Start with mass $n$ at the origin and choose $\gamma_n$ such that $\mathcal{K} \gamma_n (x,t) = 1$ if $(x,t) \neq (0,0)$, and $\mathcal{K} \gamma_n (0,0) = 1-n$. Then the odometer function $u_n$ is given by $u_n = s_n- \gamma_n$, where $s_n$ is the least supercaloric majorant of $\gamma_n$.
\label{supercal!}
\end{lemma}

Proof:
We want to prove the following:
\begin{eqnarray}
u_n + \gamma_n & = & \inf\{ f(x,t) | \mathcal{K}f \leq 0 \text{ and } f \geq \gamma_n \}
\end{eqnarray}

Let us first remark that $u_n + \gamma_n$ is indeed a supercaloric function such that $u_n + \gamma_n \geq \gamma_n$.

We recall the mass equation (\ref{important1}) as stated in terms of the discrete heat operator:
$$ \nu_n(x) - \sigma_n(x) = \mathcal{K}u_n (x). $$

Recall that the initial distribution of mass is $n$ at the origin, while the final distribution has values that are everywhere less than $1$. Hence, $u_n + \gamma_n$ is supercaloric everywhere, and since $u_n \geq 0$ by definition, it is also a majorant of $\gamma_n$.

Let us now prove that is it the smallest of these functions.
Let $f$ be a supercaloric majorant of $\gamma_n$. Then:
\begin{eqnarray}
\mathcal{K}\left( f - \gamma_n - u_n \right) & = & \mathcal{K}f - \mathcal{K}\left(u_n+\gamma_n \right) \nonumber 
\end{eqnarray}

On $D_n$, $\mathcal{K}\left(u_n + \gamma_n\right) = 0$, so that $ \mathcal{K}\left( f - \gamma_n - u_n \right) \leq 0$.
Outside $D_n$, $u_n = 0$, so that we have $f- \gamma_n - u_n = f-\gamma_n \geq 0$.

Hence, $f - \gamma_n - u_n$ is nonnegative outside $D_n$ and supercaloric inside $D_n$. The minimum principle yields that $f- \gamma_n - u_n \geq 0$ everywhere, and $f$ is indeed always greater than $u_n + \gamma_n$, which concludes the proof.

\textbf{Remark:} Lemma \ref{supercal!} gives us a way to find the odometer function for the unfair divisible sandpile model, provided that we have a function that satisfies the following conditions:
\begin{eqnarray}
\mathcal{K} \gamma_n (x,t) &=& 1 \text{ if $(x,t) \neq (0,0)$, and } \nonumber\\
\mathcal{K} \gamma_n (0,0) &=& 1-n. \label{conditions}
\end{eqnarray}

We introduce the random walk $S$ on $\Zd$, with the following law:
\begin{eqnarray}
\IP\left( S (t+1) - S(t) = \pm e_i \right) &=& \frac{1-p}{2(d-1)} \qquad \text{ for $i=1 \cdots d-1,$ and} \nonumber\\
\IP\left( S (t+1) - S (t) = e_d \right) &=& p. \nonumber\\ 
\end{eqnarray}
Our function $\gamma_n$ is easily written in terms of the random walk $S$. Define $g(x,y)$ as the expected number of times that $S$, started at $x$, visits $y$. Then consider the function:
\begin{eqnarray}
\gamma_n (z)  & = & z_d - \sum_{i=1}^{d-1}z_i^2 - n g(0,z). \nonumber\\
\end{eqnarray}
Since $\mathcal{K}g(0,z) = \delta_0(x)$, one can check that this function satisfies the conditions (\ref{conditions}).

In order to be able to use this function, we will need to evaluate the normalisation of Green's function, which is the object of the following lemma:
\begin{lemma}
For $z \in \Zd$, Define $g_n(0,z) = n g(0,z)$. Then the normalised version of $g_n$, namely $n^{-\frac{2}{d+1}}\phi_n(g_n)$, converges uniformly on compacts of  $\IR^{d-1} \times \IR^*_+$ towards the function $$G(x,t) = \frac{1}{p}    \left( \frac{\beta}{\pi t} \right)^{\frac{d-1}{2}} \exp \left( - \beta \frac{||x||^2}{t} \right).$$
\label{CVGreen}
\end{lemma}

To prove this lemma, our first step is to couple our random walk $S$ with the random walk $\tilde S$, whose steps are the sum of those of $S$ between two jumps in the direction of the drift. In other terms, the law of the increments of $\tilde S$ is that of the last position of simple random walk on $\mathbb{Z}^{d-A}$ killed at a geometric time of parameter $p$.

We calculate the covariance matrix of the walk $\tilde S$, so that we can estimate its position at a given time using the local central limit theorem (see \cite{spitzer1964principles}, \textbf{P9}). 

Let us denote $\mathbf{P}_k(0,z) = \IP \left( \tilde S_k = z \right)$. Then applying the local central limit theorem, we get that for all $(x,t) \in \IR^{d-1} \times \IR^*_+$,

\begin{eqnarray}
n^{\frac{d-1}{d+1}}\left( 2\pi t\right)^{\frac{d-1}{2}} \mathbf{P}_{\lfloor n^{\frac{2}{d+1}} t \rfloor} \left( 0, \lfloor n^{\frac{1}{d+1}} x \rfloor \right)  \rightarrow_{n \rightarrow \infty}   (2 \beta)^{\frac{d-1}{2}} \exp \left( - \frac{\beta ||x||^2}{t} \right),\nonumber
\end{eqnarray}
and the convergence is uniform on compacts of $\IR^{d-1} \times \IR^*_+$. Moreover, the difference between these quantities is of order $n^{-\frac{1}{d+1}}\frac{1}{\sqrt t}$, see for instance \cite{gnedenko1968limit}.

The probability we estimated here is that of the event that the first site visited on the $\lfloor n^{\frac{2}{d+1}} t \rfloor$-th layer in the drift direction by $S$ is $ \lfloor n^{\frac{1}{d+1}} x \rfloor$. Hence,
\begin{equation}
g\left(0,  \left( \lfloor n^{ \frac{1}{d+1}} x \rfloor, \lfloor  n^{\frac{2}{d+1}} t \rfloor \right) \right) = \sum_{z \in \mathbb{Z}^{d-1}} \mathbf{P}_{\lfloor n^{\frac{2}{d+1}} t \rfloor} \left( 0, \lfloor n^{\frac{1}{d+1}} x \rfloor +z \right)  g((z,0),0) \label{somme1}
\end{equation}

 Equation (\ref{somme1}) states that in order to get the actual number of visits to $\lfloor n^{\frac{1}{d+1}} x \rfloor$, we have to consider the sum of the hitting probabilities of $\lfloor n^{\frac{1}{d+1}} x \rfloor+ z$, where $z \in \mathbb{Z}^{d-1}$, multiplied the number of visits to $0$ of the drifted random walk started at point $(z,0)$. The hitting probabilities of $\lfloor n^{\frac{1}{d+1}} x \rfloor+ z$ with $z$ small enough to be in range of $0$ are asymptotically close to that of $\lfloor n^{\frac{1}{d+1}} x \rfloor$, so that we only have to evaluate the following sum:
\begin{equation}
\sum_{z \in \mathbb{Z}^{d-1}}  g((z,0),0) = \sum_{z \in \mathbb{Z}^{d-1}} g(0, (z,0)), \nonumber
\end{equation}
using the fact that the origin and $(z,0)$ both have coordinate $0$ in the drift direction. The summation yields exactly the expected time spent by the random walk $S$ on one given layer.

\subsubsection{Convergence of the obstacle function}
\begin{lemma}

The normalised obstacle $n^{-\frac{2}{d+1}} \phi_n(\gamma_n)$ converges uniformly on compacts $K$ of $\IR^{d-1} \times \IR$ towards $\gamma$.
\end{lemma}
The proof of this lemma is straightforward since the first two terms in $\gamma_n$ converge uniformly to those of $\gamma$ once normalised. The last term is the normalisation of Green's function, the convergence of which was proved in lemma \ref{CVGreen}.

Hence, our normalised obstacle function $n^{-\frac{2}{d+1}}\phi_n\left( \gamma_n \right)$ converges uniformly on compacts toward the function:
$$ \gamma(x,t) = t - ||x||^2 - \frac{1}{p}    \left( \frac{\beta}{\pi t} \right)^{\frac{d-1}{2}} \exp \left( - \beta \frac{||x||^2}{t} \right). $$

\subsubsection{Parabolic Potential Theory}
Before we go on, we will need a few lemmas of parabolic potential theory. Since $n^{-\frac{2}{d+1}}\phi_n(g_n)$ converges uniformly on compacts of  $\IR^{d-1} \times \IR^*_+$ towards the function $$G(0,x) = \frac{1}{p}    \left( \frac{\beta}{\pi t} \right)^{\frac{d-1}{2}} \exp \left( - \beta \frac{||x||^2}{t} \right),$$
we want to use this convergence to find a candidate to be the inverse of the heat operator, which we define in both the continuous and discrete setting. Let us define, for $z \in \Zd$ and $f$ a measurable function defined on $\IR^{d-1} \times \IR^*_+$ with compact support,
\begin{eqnarray}
G_n ( f )(z) &=& n \sum_{(y,r) \in \Zd, r>0} g(z,(y,r)) f \left( n^{-\frac{1}{d+1}}y,  n^{-\frac{2}{d+1}} r \right), \nonumber
\end{eqnarray}
and for $x \in \IR^{d-1} \times \IR^*_+$,
\begin{eqnarray}
G (f) (x) &=& \int_{(y,r) \in \IR^{d-1} \times \IR^*_+} G\left(x,(y,r)\right) f(y,r) dydr. \nonumber
\end{eqnarray}

The uniform convergence of the normalisation of $g_n$ towards $G$, together with our estimate, will ensure the convergence of the normalisation of $G_n(f)$ toward $G(f)$. This is the object of the following lemma:

\begin{lemma}
Let $f$ be a bounded function defined on $\IR^{d-1} \times \IR^*_+$ with compact support, then $$|n^{-\frac{2}{d+1}}\phi_n ( G_n(f) ) - G(f) | \rightarrow 0,$$
uniformly on compacts of $\IR^{d-1} \times \IR^*_+$.
\label{inverseK}
\end{lemma}

Suppose that $f$ is bounded by $M$ and supported on the compact $K$. We define, for $(y,r) \in K$ such that $\phi_n(y,r) \in \Zd$, the following error term:
$$ \alpha_n^x(y,r) = \left| n^{-\frac{2}{d+1}}\phi_n(g)\left(x,(y,r)\right)  -\int_{z
} G(x,z)dz \right|, \label{erreur}$$
where the integral is taken over points $z \in (y,r) + [0,n^{-\frac{1}{d+1}})^{d-1} \times [0,n^{-\frac{2}{d+1}})$.

Then we have:
\begin{eqnarray}
\left| n^{-\frac{2}{d+1}}\phi_n ( G_n(f) )(x) - G(f)(x) \right| & \leq&  \sum_{(y,r)} M  \alpha_n^x(y,r) \label{machin} 
\end{eqnarray}
where the sum is taken over points $(y,r) \in K$ such that $\phi_n(y,r) \in \Zd$.

Since we have a uniform bound on the difference $n^{-\frac{2}{d+1}}\phi_n(g)\left(x,(y,r) \right) - G\left(x,(y,r) \right)$ that has finite integral on $K$, and using our knowledge of the derivatives of $G$ to bound $G\left(x,(y,r) \right)- G\left(x,z \right)$ when $z$ is close to $(y,r)$, we can bound this sum uniformly on $K$, which concludes the proof.

The idea behind the introduction of $G(f)$ is that $G$ is in some sense the inverse of the heat operator. In fact, we do not need such a strong result, and we are content with the following property:
\begin{lemma}
Let $f$ be a positive measurable function defined on $\IR^{d-1} \times \IR^*_+$ with compact support, then $G(f)$ is supercaloric. \label{lemma2.8}
\end{lemma}

We omit the proof as it is straightforward and only relies on commuting integral signs and using the supercaloric property of $G(x,y)$ as a function of $x$.



We are now ready to prove the convergence of the odometer functions in our model.

\subsubsection{Convergence of the odometer function}
Let us consider $s_n$, the least supercaloric majorant of $\gamma_n$. The following lemma states a convergence result of the normalised version of $s_n$.
\begin{lemma}
The normalised version of the least supercaloric majorant of $\gamma_n$, namely  $n^{-\frac{2}{d+1}}\phi_n(s_n)$, converges uniformly to $s$, the least supercaloric majorant of $\gamma$, on compacts of $\IR^{d-1} \times \IR^*_+$.
\label{CVUu}
\end{lemma}

We begin our proof by pointing out that our definition of the least supercaloric majorant depends on the context: while $s_n$ is defined in the discrete setting, $s$ is defined in the continuous.

Let $K$ be a compact of $\IR^{d-1} \times \IR^*_+$, and $K_j$ an increasing sequence of compacts of $\IR^{d-1} \times \IR^*_+$ such that $\cup_j K_j = \IR^{d-1} \times \IR^*_+$, and $K_j \subset \mathring{K_{j+1}}$. We also define their discrete counterparts $K_j^n$ as the set of points $(x,t)$ such that $\left(n^{-\frac{1}{d+1}} x, n^{-\frac{2}{d+1}} t\right) \in K_j$.

We define the following quantities:
\begin{eqnarray}
s_n^{K_j} &=& \inf\{ f(x,t) | \mathcal{K}f \leq 0 \text{ globally, and } f \geq \gamma_n  \text{ on } K_j^n \}, \nonumber\\
s^{K_j} &=& \inf\{h(x,t)|  h(x,t)\text{ is globally supercaloric and } h\geq \gamma \text{ on } K_j\}. \nonumber
\end{eqnarray}

Remark that both functions are increasing in $j$, so it is a consequence of Dini's theorem that $s^{K_j}$ and $\phi_n(s_n^{K_j})$ converge uniformly on $K$ towards $s$ and $\phi_n(s_n)$, respectively, as $j$ tends to infinity.

We now write that:
$$|n^{-\frac{2}{d+1}}\phi_n(s_n) - s| \leq |n^{-\frac{2}{d+1}}\phi_n(s_n) - n^{-\frac{2}{d+1}}\phi_n(s_n^{K_j})| +|n^{-\frac{2}{d+1}}\phi_n(s^{K_j}_n)-s^{K_j}|+|s^{K_j}-s|$$

Set $j$ big enough so that $K \subset K_j$, and $n^{-\frac{2}{d+1}}|\phi_n(s_n) - \phi_n(s_n^{K_j})| +|s^{K_j}-s| \leq \epsilon$ on $K$. It now remains to show that $n^{-\frac{2}{d+1}}\phi_n(s^{K_j}_n)$ converges uniformly on $K$ towards $s^{K_j}$ as $n$ tends to infinity.

We will now proceed to smooth the function $s^{K_j}$ in order to show that it is close to a function $\tilde s^{K_j}$ such that $\mathcal{K} \tilde s^{K_j}$ is small. 

We first remark that $s^{K_j}$ is continuous on $K$. Indeed, as it is an infimum of continuous functions, it is also upper semi-continuous. Moreover, if we define $\omega(\gamma, r)$ as the continuity modulus of $\gamma$ on $K_{j+1}$, we have, on $K_j$,
\begin{equation}
\mathfrak{A}_r \left(s^{K_j} \right) \geq \mathfrak{A}_r \left( \gamma \right) \geq \gamma - \omega(\gamma, r),
\end{equation} 
so that $\mathfrak{A}_r \left(s^{K_j} \right) + \omega(\gamma, r)$ is continuous, supercaloric, and lies above $\gamma$ on $K_j$. Hence, $$\mathfrak{A}_r \left(s^{K_j} \right) \leq s^{K_j} \leq \mathfrak{A}_r \left(s^{K_j} \right) + \omega(\gamma, r).$$

Since $\gamma$ is uniformly continuous on $K_{j+1}$, $\mathfrak{A}_r \left(s^{K_j}  \right) \rightarrow s^{K_j}$ as $r \rightarrow 0$. Since $s^{K_j}$ is supercaloric, this is an increasing limit. As an increasing limit of continuous functions, $s^{K_j}$ is also lower semi-continuous.

We now define, like in \cite{Evans}, Appendix C, for $x \in K$, $\lambda >0$,
$$ \tilde s^{K_j} (x) = \int_{\IR^d} s^{K_j}(y) \lambda^{-d} \eta \left( \frac{x-y}{\lambda} \right), $$
where $\eta$ is defined as follows:
$$ \eta(x) = \left\{
          \begin{array}{ll}
            C \exp\left( \frac{1}{|x|^2-1} \right), & \qquad |x| <1 \\
            0, & \qquad |x| \geq 1, \\
          \end{array}
        \right.$$
with $C$ such that $\int_{\IR^d} \eta$=1 and $\lambda$ is taken sufficiently small to ensure the definition of the values of $s^{K_j}$ involved. Then $\tilde s^{K_j}$ is still supercaloric, and if $\lambda$ is small enough, it satisfies $ |\tilde s^{K_j} - s^{K_j}| \leq \epsilon$ on $K$.

Let us define the following discrete function:
$$(\forall k,l) \in \mathbb{Z}^{d-1} \times \mathbb{Z}, \quad q_n(x,t) = \tilde s^{K_j} \left( n^{-\frac{1}{d+1}}x, n^{-\frac{2}{d+1}} t\right) - n^{-\frac{5}{2(d+1)}}|x|^2.$$
Since $\tilde s^{K_j}$ is smooth and supercaloric, it satisfies the partial differential inequation:
$$\mathfrak{K}(s^{K_j}) \leq 0.$$ 
Moreover, since $\tilde s^{K_j}$ is smooth, Taylor's formula yields that, on $K$,  
$$\left\| n^{\frac{2}{d+1}} \mathcal{K}\left( s^{K_j}\left(n^{-\frac{1}{d+1}}x, n^{-\frac{2}{d+1}} t\right) \right) - \mathfrak{K}\left( \tilde s^{K_j} \right) \right\| \leq A n^{\frac{3}{d+1}},$$
where $A$ is chosen to be the maximum of the norms of third derivatives in space dimensions, and second derivative in time of $\tilde s^{K_j}$.

Hence, the smoothness and supercaloric property of $\tilde s^{K_j}$ ensure that, for $n$ large enough:
$$ \mathcal{K} q_n(x,t) \leq 0.$$

When $n$ is large enough, $\phi_n(q_n)$ is close to $\tilde s^{K_j}$, which is in turn close to $s^{K_j}$, which is greater than $\gamma$, itself close to $n^{-\frac{2}{d+1}}\phi_n(\gamma_n)$ on $K_j$. To sum up, the following inequalities hold on $K_j$:

$$ \phi_n(q_n) > \tilde s^{K_j} - \epsilon > s^{K_j} - 2\epsilon > \gamma -2\epsilon > n^{-\frac{2}{d+1}}\phi_n(\gamma_n)-3\epsilon. $$

Hence $\phi_n(q_n) +3\epsilon$ is a supercaloric majorant of $n^{-\frac{2}{d+1}}\phi_n(\gamma_n)$ on $K_j$, so $\phi_n(s_n^{K_j}) \leq \phi_n(q_n) + 3\epsilon \leq \tilde s^{K_j} + 3\epsilon \leq s^{K_j} + 4\epsilon$.

We will now prove the converse inequality.

Define the following functions:
\begin{eqnarray}
h_j^n &=& - \mathcal{K}\left( s^{K_j}_n \mathbf{1}_{K_{j+1}^n} \right), \nonumber\\
G (h_j^n) (x) &=& \int_{(y,r) \in \IR^{d-1} \times \IR^*_+} G\left(x,(y,r)\right) h_j^n(\lfloor n^{\frac{1}{d+1}}y \rfloor, \lfloor n^{\frac{2}{d+1}}r \rfloor)dydr. \nonumber
\end{eqnarray}

We will start by proving that on $K_{j+1}^n$,
$$\left| \mathcal{K} s^{K_j}_n \right| \leq 1. \label{petith}$$
To see this, a discrete reasoning is necessary. Indeed, let us look at one particular point $z$, and suppose $\mathcal{K} s^{K_j}_n(z) < -1$. Consider the function $f$ that coincides with $ s^{K_j}_n $ on every point but $z$, and adjust the value $f(z)$ such that $\mathcal{K} f(z)=-1$. Remark that this implies $f(z) < s^{K_j}_n(z)$. We will now prove that $f$ is a supercaloric majorant of $\gamma_n$. 
First, it is supercaloric at point $z$ because $\mathcal{K} f(z)=-1$ by definition, and on other points, since $f(z) \leq s^{K_j}_n(z)$ implies that for $y \neq z$, $\mathcal{K} f(y) \leq \mathcal{K} s^{K_j}_n(z)$. Note that, as an infimum of discrete supercaloric functions, $s^{K_j}_n(z)$ is itself supercaloric, so that for $y \neq z$, $\mathcal{K} f(y) \leq \mathcal{K} s^{K_j}_n(z) \leq 0$.
The function $f$ is also a majorant of $\gamma_n$, which is true at points $y \neq z$ by definition. At point $z$, we know that $\mathcal{K} f = -1 = \mathcal{K} \gamma_n$, so that the minimum principle applied locally between $z$ and its neighbors to the function $f - \gamma_n$ guaranties $f(z)-\gamma_n(z) \geq 0$.
Hence, $f$ is a supercaloric majorant of $\gamma_n$ on $K_{j+1}^n$, so that $s_n^{K_j} \leq f$, which is a contradiction. Hence, n $K_j^n$,
$$\left| \mathcal{K} s^{K_j}_n \right| \leq 1.$$

On the other hand, since $-G_n$ inverses $\mathcal{K}$ exactly, $n^{-\frac{2}{d+1}}s^{K_j}_n = G_n(h_j^n)$ on all non-boundary points of $K_{j+1}^n$.

We want to argue that $G(h_j^n)(x)$ is a good approximation of $n^{-\frac{2}{d+1}} \phi_n(s^{K_j}_n)(x)$. We separate the error as follows, for $x \in K_j^n$,
\begin{eqnarray}
\left| G(h_j^n)(x) - n^{-\frac{2}{d+1}} \phi_n(s^{K_j}_n)(x) \right| \leq A+B \nonumber
\end{eqnarray}
where 
\begin{eqnarray}
A &\leq & \sum_{(y,r) \in K_{j+1}^n}  |h_j^n(y,r)| \alpha_n^x(y,r) \nonumber\\
B &\leq & B_0 n^{-\frac{2}{d+1}} \sum_{(y,r), (z,t)} s_n^{K_j}(y,r) \alpha_n^x(z,t) - s_n^{K_j}(z,t) \alpha_n^x(y,r), \nonumber
\end{eqnarray}
where the second sum is taken over points $(y,r) \sim (z,t)$ such that $(y,t) \in K^n_{j+1}$, but $ (z,t) \notin K^n_{j+1}$, and $B_0$ is a suitable constant depending only on the dimension. The term $B$ estimates the part of the error that stems from using $s^{K_j}_n \mathbf{1}_{K_{j+1}^n}$ rather than $s^{K_j}_n$ in the definition of $h^n_j$.
To bound $A$, we use equation (\ref{petith}):
\begin{equation}
A \leq  \sum_{(y,r) \in K_{j+1}^n}  \alpha_n^x(y,r) \nonumber
\end{equation}
We then find the same sum as in the proof of lemma \ref{inverseK}.

Bounding $B$ is a little trickier. We proceed in the following way:
\begin{eqnarray}
\left| s_n^{K_j}(y,r) \alpha_n^x(z,t) - s_n^{K_j}(z,t)\alpha_n^x(y,r) \right| &\leq&  \left| s_n^{K_j}(y,r) \right| \left|\alpha_n^x(z,t) - \alpha_n^x(y,r) \right| \nonumber\\
 & & + \left|\alpha_n^x(y,r)  \right| \left| s_n^{K_j}(y,r) - s_n^{K_j}(z,t) \right|. \nonumber
\end{eqnarray}

 We will rely on three arguments: first, remark that both $(y,r)$ and $(z,t)$ are safely away from $x$, because of our condition $K_j \subset \mathring{K_{j+1}}$ with $(y,r)$ and $(z,t)$ neighbors, a short computation shows that:
$$\alpha_n^x(z,t) - \alpha_n^x(y,r) = o\left(n^{-\frac{2}{d+1}} \right).$$

Secondly, it follows from the definition of $s_n^{K_j}$ and the convergence of $n^{-\frac{2}{d+1}}\phi_n(\gamma_n)$ that on $K_j^n$,
$$ s_n^{K_j}(z) \leq C n^{\frac{2}{d+1}}, $$
where $C$ is a constant depending on $K_j$.

Last, we use the fact that $\gamma$ is uniformly continuous on $K_j$ to argue that the increments of $\gamma_n $ between two neighboring points are at most of order $ \epsilon n^{\frac{2}{d+1}}$ for $n$ large enough, and this property is in turn transferred to $s_n^{K_j}$ (since $s_n^{K_j}(. + e_i)+\epsilon $ is a supercaloric majorant of $\gamma_n$).

Hence, our terms are bounded in a satisfactory way, so that on $K_j$, $G(h_n)$ is uniformly close to $n^{-\frac{2}{d+1}}\phi_n(s^{K_j}_n)$, which is bigger than $n^{-\frac{2}{d+1}}\phi_n(\gamma_n)$, which is uniformly close to $\gamma$. To sum up, on $K_j$,

$$  G(h_n)  > n^{-\frac{2}{d+1}}\phi_n(s^{K_j}_n) - \epsilon > n^{-\frac{2}{d+1}}\phi_n(\gamma_n) - \epsilon > \gamma - 2 \epsilon. $$ 

Moreover, lemma \ref{lemma2.8} states that $G(h_n)$ is supercaloric everywhere, so that $G(h_n) + 2 \epsilon \geq s$ and we can conclude that $n^{-\frac{2}{d+1}}s^{K_j}_n+3\epsilon \geq s^{K_j}$.

The convergence of $n^{-\frac{2}{d+1}}\phi_n(s^{K_j}_n)$ towards $s^{K_j}$ is thus uniform on $K$. Recall that the  convergences of $s^{K_j}$ towards $s$ and $s_n^{K_j}$ towards $s_n^{K_j}$ are guaranteed to be uniform by Dini's theorem, which concludes the proof.

\begin{lemma}
The normalised odometer function $n^{-\frac{2}{d+1}} \phi_n(u_n)$ converges towards $u=s- \gamma $ uniformly on compacts of $\IR^{d-1} \times \IR^*_+$.

Moreover, the function $u$ and the family $n^{-\frac{2}{d+1}} \phi_n(u_n)$ have bounded (respectively uniformly bounded) integrals over any compact of $\IR^{d-1} \times \IR_+$.
\label{helper}
\end{lemma}

The first statement is a consequence of the convergence properties of $s_n$ and $\gamma_n$, while the second relies on the fact that $t- ||x||^2$ is a supercaloric majorant of $\gamma$ (resp. $\gamma_n$). Therefore, for $x \in \Zd$,
$$ s_n(x) - \gamma_n(x) \leq n g(0,x), $$
and for $(x,t) \in \IR^{d-1} \times \IR^*_+$,
$$s - \gamma \leq \frac{1}{p} \left(\frac{\beta}{\pi t} \right)^{\frac{d-1}{2}} \exp\left( - \frac{\beta ||x||^2}{t} \right).$$

A quick computation finishes the proof.

\subsubsection{Convergence of Domains}
\begin{lemma}
Let $K$ be a compact of $\IR^{d-1} \times \IR^*_+ $. Then, on $K$, the normalised unfair divisible sandpile cluster $\phi_n\left(D_{n}\right)$ converges to $D$ with respect to the Hausdorff distance.
\end{lemma}

Set $\epsilon >0$, and define $D_{\epsilon}$ as the inside $\epsilon$-neighbourhood of $D$, that is to say,
$$ D_{\epsilon} = \{ x \in D, \B_{x, \epsilon} \subset D \}.$$

Since $D$ is the non-coincidence set for the caloric obstacle problem with obstacle $\gamma$, $s-\gamma$ is strictly positive on $\bar{D_\epsilon}$. Since $s-\gamma$ is uniformly continuous on $K \cap \bar{D_{\epsilon}}$, let us define $$ \beta = \min_{z \in K \cap \bar{D_\epsilon}} ( s(z) - \gamma(z) ) >0.$$

Then, for $n$ large enough, the uniform convergences of the normalised obstacle and majorant guarantee that $|\gamma - n^{-\frac{2}{d+1}}\phi_n(\gamma_n)| < \frac{\beta}{2}$, and $|s -n^{-\frac{2}{d+1}}\phi_n(s_n)| < \frac{\beta}{2}$ on $K \cap \bar{D}$. This yields the following :
$$ n^{-\frac{2}{d+1}}\phi_n(s_n) - n^{-\frac{2}{d+1}}\phi_n(\gamma_n) > s - \gamma - \beta > 0 \text{ on } K \cap D_{\epsilon}, $$
so that $K \cap D_{\epsilon} \subset K\cap \phi_n(D_{n})$ for all but finitely many $n$.

Let $D^{\epsilon}$ be the outside $\epsilon$-neighbourhood of $D$.
Let $(x_n, t_n) \in D_{n}$, and suppose that $(n^{-\frac{1}{d+1}}x_n,n^{-\frac{2}{d+1}} t_n)$ converges to $(x_0,t_0) \in \IR^{d-1} \times \IR^*_+$. It is sufficient to prove that $(x_0,t_0) \in D^{\epsilon}$.

Define the discrete cylinder segment $\mathcal{C} = \mathbf{B}\left( x_n, n^{\frac{1}{d+1}}\frac{\epsilon}{2} \right)  \times [\lfloor t_n - n^{\frac{2}{d+1}}\frac{\epsilon^2}{4} \rfloor, \lfloor t_n \rfloor]$.

On $\mathcal{C} \cap D_{n}$, define the following function :
$$w_n(x,t) = u_n(x,t) - \left( |x-x_n|^2 - (t-t_n) \right)$$

Then $\mathcal{K}(w_n) = \mathcal{K}u_n -1 \geq 0$, so that $-w_n$ is discrete supercaloric.

Hence, $w_n$ satisfies the parabolic maximum principle and realises its maximum on the parabolic boundary of $\mathcal{C} \cap D_{n}$.

Since $w_n(x_0, t_0) = u_n(x_0,t_0) >0$, and $u_n = 0$ on $\partial D_{n}$, the maximum cannot be taken on $\partial D_{n}$.

Hence $w_n$ takes its maximum on the parabolic boundary of $\mathcal{C}$, $\partial_p \mathcal{C}$. This set can be decomposed as follows:
$$\partial_p \mathcal{C} = \Lambda_1 \cup \Lambda_2,$$
where 
\begin{eqnarray}
\Lambda_1 & = & \mathbf{B}\left( x_n, n^{\frac{1}{d+1}}\frac{\epsilon}{2} \right) \times \{ \lfloor t_n - n^{\frac{2}{d+1}}\frac{\epsilon^2}{4} \rfloor \} \text{, and } \nonumber\\
\Lambda_2 & = & \partial \mathbf{B}\left( x_n, n^{\frac{1}{d+1}}\frac{\epsilon}{2} \right) \times \left[\lfloor t_n- n^{\frac{2}{d+1}}\frac{\epsilon^2}{4} \rfloor, \lfloor t_n \rfloor \right] \nonumber
\end{eqnarray}

Now, on $\Lambda_1,$ we have:
$$ w_n(x,t) \leq u_n(x,t) - n^{\frac{2}{d+1}}\frac{\epsilon^2}{4}, $$
and, similarly, on $\Lambda_2,$
$$ w_n(x,t) \leq u_n(x,t) - n^{\frac{2}{d+1}}\frac{\epsilon^2}{4}.$$

Let $(y_n, r_n)$ be the point where $w_n$ takes its maximum. Then,
\begin{eqnarray}
u_n(x_n,t_n) = w_n(x_n,t_n) & \leq & w_n(y_n,r_n) \nonumber\\
                            & \leq & u_n(y_n,r_n) \nonumber
\end{eqnarray}
This gives the following:
$$ u_n(y_n,r_n) \geq  u_n(x_n, t_n) + n^{\frac{2}{d+1}}\frac{\epsilon^2}{4}.$$
Let us extract a subsequence of $(n^{-\frac{1}{d+1}}y_n,n^{-\frac{1}{d+1}}r_n)$ that converges to a point $(y,r)$. Since $n^{-\frac{1}{d+1}} \phi_n(u_n)$ converges to $u$ uniformly on compacts, the limit $u(y,r)$ is strictly positive, so that $(y,r) \in D $.

However, point $(y,r)$ is within distance $\epsilon$ of $(x_0, t_0)$, so we can conclude that $(x_0, t_0) \in D^{\epsilon}$, and $D_{n} \subset D^{\epsilon}$.

\subsection{Regularity of $D$}

\subsubsection{Differential Equation approach}

So far, we have considered $u$ only as the limit of the discrete odometer, but it will be useful to see that it is also the solution of the continuous version of the equation that defines $u_n$. This is the object of the following lemma:

\begin{theorem}
The normalised limit $u$ of the odometers solves the following partial differential equation in the weak (distributional) sense:
\begin{equation}
 \frac{1-p}{2(d-1)}\Delta u - p \frac{\partial u}{\partial t} = \mathbf{1}_{u>0} -  \delta_0, \label{PDE1}
\end{equation}
where $\delta_0$ is the Dirac measure at the origin.
\end{theorem}

Consider a discrete function $\eta$ that has compact support on $\Zd$. Since $u_n$ solves the discrete equation $\mathcal{K} u_n =  \nu_n -  n\delta_0, $ we have:
\begin{eqnarray}
\sum_{z \in \Zd} \eta (z)  \mathcal{K} u_n (z)  & = & \sum_{z \in \Zd} \eta(z)  \nu_n(z) -  n \eta(0,0) \label{discretePDE1}
\end{eqnarray}
When changing the variables so as to report the operation on $\eta$, only the sign of the drift coordinate is modified, so that we define the discrete operator $\tilde{\mathcal{K}}$ as the operator $\mathcal{K}$ with reversed time:
$$  \tilde{\mathcal{K}}f(z) = \frac{1-p}{2(d-1)} \sum_{||y-z|| =1, (y-z)\perp e_d} \left( f(y) - f(z) \right) + p \left( f( z  ) - f(z- e_d) \right). $$

Equation (\ref{discretePDE1}) can now be written in the following format:
\begin{eqnarray}
\sum_{z \in \Zd} u_n (z)  \tilde{\mathcal{K}}\eta(z)  & = & \sum_{z \in \Zd} \eta(z)  \nu_n(z) -  n \eta(0,0). \label{discretePDE2}
\end{eqnarray}

Let us now consider a test function $h: \IR^{d-1} \times \IR \rightarrow \IR$ that is $C^{\infty}$ and has compact support. Define its discrete counterpart $h_n: \Zd \rightarrow \IR$ as:
$$ h_n(z_1,...,z_d) = h\left( \left( n^{-\frac{1}{d+1}} z_i \right)_{{i} \in \{ 1,...,d-1\}}, n^{-\frac{2}{d+1}} z_d \right).$$
Then applying equation (\ref{discretePDE2}) yields:
\begin{eqnarray}
\frac{1}{n}\sum_{z \in \Zd} u_n (z)  \tilde{\mathcal{K}}h_n(z)  & = & \frac{1}{n}\sum_{z \in \Zd} h_n(z)  \nu_n(z) -  h(0,0). \nonumber\\
\frac{1}{n}\sum_{z \in \Zd} \left( n^{-\frac{2}{d+1}} u_n (z) \right)  \left( n^{\frac{2}{d+1}}\tilde{\mathcal{K}}h_n(z) \right)  & = & \frac{1}{n}\sum_{z \in \Zd} h_n(z) \nu_n(z)  -  h(0,0). \label{discretePDE3}
\end{eqnarray}
Note that $\nu_n$ is equal to $1$ on $D_n$, and has values between $0$ and $1$ at distance $1$ from $D_n$. At distance $2$ or more from $D_n$, $\nu_n = 0$. Note that Theorem \ref{sandpile} proves that $\phi_n\left(D_n\right)$ converges to $D$ on compacts $K \subset \IR^{d-1} \times \IR^{*}_+$, and lemma \ref{helper} proves that $n^{-\frac{2}{d+1}} \phi_n\left(  u_n \right)$ converges to $u$ uniformly on $K$. Hence, letting $n$ go to infinity in (\ref{discretePDE3}), the regularity of $h$ and its derivatives, together with the bounded integral property stated in lemma \ref{helper}, ensures that we can apply the dominated convergence theorem, which yields:
$$\int_{\IR^{d-1} \times \IR^{*}_+} u \left( \frac{(1-p)}{2(d-1)} \Delta h + p\frac{\partial h}{\partial t} \right) = \int_{D} h - h(0,0), $$
and equation (\ref{PDE1}) holds in the distributional sense.
\newline

This characterisation of $u$ as a solution to a parabolic free boundary problem yields numerous analytical results on the function $u$ and the set $D=\{ u>0 \}$. The following proposition states the most important ones for our problem:
\begin{proposition}
Let $\phi$ be a $C^{\infty}$ function of time and space such that
$$ \frac{1-p}{2(d-1)} \Delta \phi + p \frac{\partial \phi}{ \partial t} =0. $$
Then the following mean value property holds:
$$ \int_{D} \phi(z,t) dzdt = |D| \phi(0). $$
Moreover, the limit $u$ of the odometer is continuous in both time and space coordinates, and the set $\partial D$ has (d-dimensional) Lebesgue measure $0$. 
\label{regularity}
\end{proposition}

The first part of the proposition is a direct consequence our characterisation of $u$ as a solution to the equation (\ref{PDE1}) in the distributional sense. The second part of the property states regularity results on such a solution. For the proof of this powerful result in parabolic potential theory, we refer the reader to \cite{caffarelli2004regularity}. For an in-depth study of parabolic free boundary problems and their regularity, see the detailed book by Friedman \cite{friedman1982variational}.


\section{Drifted iDLA}
In this section, we prove the almost sure convergence of the normalised drifted iDLA cluster towards $D$, the limiting shape of the unfair divisible sandpile model. For reasons of simplicity and readability, we will start by using a random walk whose restriction to the direction of the drift is an increasing process.
\subsection{The model}

Let $ (S^j)_{j \in \mathbb{N}}$ a sequence of independent random walks on $\Zd$, with the following law:
\begin{eqnarray}
\IP\left( S^j (t+1) - S^j(t) = \pm e_i \right) &=& \frac{1-p}{2(d-1)} \qquad \text{ for $i=1 \cdots d-1,$ and} \nonumber\\
\IP\left( S^j (t+1) - S^j (t) = e_d \right) &=& p. \nonumber\\ 
\end{eqnarray}

The iDLA cluster is built recursively using the random walk $S^j$ in the following fashion.
We start with an empty set. At step $j$, suppose that we have built $A(j-1)$. We launch the random walk $S^j$ at the origin, which we let evolve until it exits $A(j-1)$. The first site visited outside the cluster is then added to the cluster.

Formally, let us define the cluster $A(n)$ and stopping times $(\sigma_k)_{k \in \mathbb{N}}$ recursively in the following way:
\begin{eqnarray}
\nu_1 & = & 0 \text{,}\nonumber\\
A(1)     & = & \{ 0 \} = \{ S^1(\sigma_1) \} \text{, and }\nonumber\\
\forall j >1, \text{ } \nu_j & = & \inf \{ t \geq 0 : S^j(t) \not \in A(j-1) \} \text{,}\nonumber\\
                          A(j)  & = & A(j-1) \cup \{ S^j(\sigma_j) \} \text{.}\nonumber
\end{eqnarray}

A few simple facts can be stated about the cluster. It is a random increasing set that contains the origin. We have $\sharp A(j) = j$, and given the form of the law of the random walks, the intersection of the cluster with $\mathbb{Z}^{d-1} \times \left( - \mathbb{N}^* \right)$ is empty.

\subsection{Limiting shape}

We first prove internal convergence of our cluster intersected with compacts of $\IR^{d-1} \times \IR^*_+$. We then use this convergence, together with the fact that it yields control over a proportion of the particles arbitrarily close to $1$, to prove external convergence as well.

\begin{theorem}
[Internal convergence towards $D$, on compacts]

Let $ A(n)$ be the drifted iDLA cluster, and $K$ a compact subset of $\IR^{d-1} \times \IR^*_+$.

Then, for every $\epsilon >0$, 
$$ D_{\epsilon} \cap K  \subset \phi_n( A(n) ) \cap K,$$
for all but finitely many $n$, almost surely.

\label{idla-int}
\end{theorem}

Recall that we build the cluster using the family of independent random walks $(S^i)_{i \in \{1 \cdots n\}}$. For $z \in \Zd$, we define the following stopping times:
\begin{eqnarray}
\tau^i_z         & = &\inf\{ t \geq 0, S^i(t) = z \}, \nonumber\\
\tau^i_{D_n}& = &\inf\{ t \geq 0, S^i(t),  \notin D_n \}, \nonumber\\
\nu^i            & = &\inf\{ t \geq 0, S^i(t) \notin A(i-1) \}. \nonumber
\end{eqnarray}
We introduce the following counting random variables :
\begin{eqnarray}
\mathcal{L}_{ n}(z) &=& \sum_{i=1}^{ n} \mathbf{1}_{\nu^i < \tau^i_z < \tau^i_{D_n}}, \nonumber\\
\mathcal{M}_{ n}(z) &=& \sum_{i=1}^{ n} \mathbf{1}_{\tau^i_z < \tau^i_{D_n}}, \nonumber\\
\mathcal{N}_{ n}(z) &=& \sum_{i=1}^{ n} \mathbf{1}_{\tau^i_z \leq \tau^i_{D_n} \wedge \nu^i}. \nonumber
\end{eqnarray}

Since $\mathcal{N}_{ n}^z$ measures the number of particles that pass through $z$ before either adding to the cluster or leaving $D_n$, if $\mathcal{N}_{ n}(z) >0$, then point $z$ is in the cluster $A(n)$.

Remark that $\mathcal{N}_{ n}(z) \geq \mathcal{M}_{ n}(z) - \mathcal{L}_{ n}(z)$.

The strategy of the proof will be to study these random variables. After showing that their expected value is different enough, and that they are close enough to their respective expected values, we will deduce that their difference is suitably big as $n$ tends to $\infty$.

First, we use the Markov property of our random walk to argue that summing over all walks that add to the cluster before exiting $D_n$ can be re-written as a sum over points of $D_n$. Let us construct the following family of random walks $S^y$, for $y \in D_n$ as follows. If there is an index $i \in \{1 \cdots n \}$ such that $y = S^i(\nu_i)$, then we define $S^y(t) = S^i(t+\nu_i)$. Remark that the aggregation rule guarantees that there is at most one such index. If there is no such index, then we define $S^y$ as a random walk started at $y$, with the same increments as $S^1$ for instance, and independent of all the other random walks.

We define the corresponding stopping times:

\begin{eqnarray}
\tau^y_z         & = &\inf\{ t \geq 0, S^y(t) = z \}, \nonumber\\
\tau^y_{D_n}& = &\inf\{ t \geq 0, S^y(t),  \notin D_n \}. \nonumber\\
\end{eqnarray}

Then we have the following inequality:
\begin{eqnarray}
\mathcal{L}_{ n}(z) &=& \sum_{i=1}^{ n} \mathbf{1}_{\nu^i < \tau^i_z < \tau^i_{D_n}} \nonumber\\
 & \leq & \sum_{y\in D_n}\mathbf{1}_{\tau^y_z < \tau^y_{D_n}} = \tilde{ \mathcal{L}}_{ n}(z), \nonumber 
\end{eqnarray}
which holds because every term in the first sum can be found once in the second sum, with possible additional terms. Remark that $\tilde{\mathcal{L}}_n(z)$ is now a sum of independent indicator random variables.

Let us compute the expectations of these variables:
\begin{eqnarray}
\IE\left( \mathcal{M}_{ n}(z) \right) & = & n \frac{g_{n,D_n}(0,z)}{g_{n,D_n}(z,z)} \nonumber\\
\IE\left( \tilde{\mathcal{L}}_{ n}(z) \right) & = & \frac{1}{g_{n,D_n}(z,z)} \sum_{y\in D_n} g_{n,D_n}(y,z), \nonumber
\end{eqnarray}
where $g_{n,D_n}$ is the Green function of a random walk $S$ stopped when it exits $D_n$.

Define  $S^{\leftarrow}$ a random walk with opposite drift, and $g_{n,D_n}^{\leftarrow}$ the corresponding stopped Green function. For any pair of points $y$ and $z$, counting the trajectories from $y$ to $z$ and from $z$ to $y$ yields that $g_{n,D_n}(y,z)= g_{n,D_n}^{\leftarrow}(z,y)$. Hence,
\begin{eqnarray}
\IE\left( \tilde{\mathcal{L}}_{ n}(z) \right) & = & \frac{1}{g_{n,D_n}(z,z)} \sum_{y\in D_n} g^{\leftarrow}_{n,D_n}(z,y) \nonumber\\
&=& \frac{\IE(\tau^{\leftarrow}_{z,n,D_n})}{g_{n,D_n}(z,z)}, \label{oppositedrift}
\end{eqnarray}
where $\tau^{\leftarrow}_{z,n,D_n}$ is the time at which the random walk $S^{\leftarrow}$, started at $z$ exits $D_n$ (i.e. the total time spent by the walk in $D_n$).

Define $f_{n, D_n}(z) = g_{n, D_n}(z,z) \IE\left( \mathcal{M}_{ n}(z) - \tilde{\mathcal{L}}_{ n}(z) \right).$ Then,
\begin{eqnarray}
f_{n,D_n}(z) & = & g_{n, D_n}(z,z) \left( \sum_{i=1}^{ n} \IP\left( \tau^i_z < \tau^i_{D_n} \right) - \sum_{y\in D_n} \IP\left(\tau^y_z < \tau^y_{D_n}\right)\right) \nonumber\\
 & = & g_{n, D_n}(z,z) \sum_{y\in D_n}( \delta_0(y)  n -1) \IP ( \tau^y_z < \tau^y_{D_n}). \nonumber\\
 & = &  \sum_{y\in D_n}( \delta_0(y)  n -1)g_{n, D_n}(y,z) \nonumber
\end{eqnarray}

Hence, $f_{n, D_n}$ is a solution to the discrete parabolic PDE :
$$\begin{array}{rlll}
\mathcal{K} f_{n,D_n}(z) & = &  1 - n\delta_0(z) & \text{ $\forall z \in D_n$,} \\
 f_{n,D_n}(z) & = & 0 & \text{ $\forall z \notin D_n$.} 
\end{array} $$

Recall that $u_n$, the odometer function for the unfair divisible sandpile, satisfies the same equation, on the interior of $D_n$. This means that $f_{n, D_n} - u_n$ is a discrete caloric function on the interior of $D_n$. Moreover, for $z$ on the interior boundary of $D_n$, since $z$ has at least one neighbor who never toppled, it means that the total mass emitted from $z$ towards this neighbor is less than $1$, which means we have:
$$ u_n(z) \leq 2d. $$

 Hence, for all $z \in \Zd$,
$$ f_{n,D_n}(z) \geq u_n(z) -2d$$




Let us define the following minimal value of $u$:
$$ \beta = \inf_{x \in D_{\epsilon} \cap K} u(x). $$

Since $u$ is continuous and $K$ is a compact, the infimum is actually a minimum, and we have $\beta >0$.

Since $n^{-\frac{2}{d+1}} \phi_n(u_n)$ converges uniformly to $u$ on $K$, we can choose $n$ large enough so that for all points $z \in \Zd$ such that:
$$\left( \left(n^{-\frac{1}{d+1}}z_i\right)_{i\in \{1 \cdots d-1\}}, n^{-\frac{2}{d+1}} z_d \right) \in D_{\epsilon} \cap K,$$
we have, for $n$ large enough,
$$u_n(z) \geq n^{\frac{2}{d+1}}\frac{\beta}{2} $$


This means the following relation on the expected values of $\tilde{\mathcal{L}}_{ n}(z) $ and $\mathcal{M}_{ n}(z) $ holds for $n$ large enough:
\begin{eqnarray}
\IE\left( \mathcal{M}_{ n}(z) - \tilde{\mathcal{L}}_{ n}(z) \right) & \geq & \frac{1}{2} \frac{\beta n^{\frac{2}{d+1}}}{g_{n,D_n}(z,z)}
\end{eqnarray}

Using equation (\ref{oppositedrift}), we get, for $n$ large enough:
\begin{eqnarray}
\IE\left(\mathcal{M}_{ n}(z) \right) & \geq & \left( 1 + \frac{\beta n^{\frac{2}{d+1}}}{2  \IE \left( \tau^{\leftarrow}_{z,n,D_n}\right)} \right) \IE \left( \tilde{\mathcal{L}}_{ n}(z) \right) \label{majoration}
\end{eqnarray}

Since$\left( \left(n^{-\frac{1}{d+1}}z_i\right)_{i\in \{1 \cdots d-1\}}, n^{-\frac{2}{d+1}} z_d \right) \in D_{\epsilon} \cap K,$ we have for a suitably large constant $C'$, uniformly in $z$,
$$ \IE \left( \tau^{\leftarrow}_{z,n,D_n}\right) \leq C' n^{\frac{2}{d+1}}.$$

Hence, equation (\ref{majoration}) can be written as
\begin{eqnarray}
\IE\left(\mathcal{M}_{ n}(z) \right) & \geq & \left( 1+ \kappa \right) \IE \left( \tilde{\mathcal{L}}_{ n}(z) \right), \nonumber
\end{eqnarray}
where $\kappa$ is a strictly positive constant.

It follows that we can choose $\lambda >0$ such that $$(1+ \lambda) \IE \left( \tilde{\mathcal{L}}_{ n}(z) \right) < (1-\lambda) \IE\left(\mathcal{M}_{ n}(z) \right),$$ and apply the following large deviation principle, which holds for sums of independent indicator variables :
\begin{eqnarray}
\IP\left( \mathcal{L} \geq (1+\lambda) \IE\left( \mathcal{L} \right)\right) &<& 2 e^{-c_{\lambda} \IE\left( \mathcal{L} \right) } \nonumber\\
\IP\left( \mathcal{M} \geq (1-\lambda) \IE\left( \mathcal{M} \right)\right) &<& 2 e^{-c_{\lambda} \IE\left( \mathcal{M} \right) } \nonumber
\end{eqnarray}
(See \cite{alon2008probabilistic}, Cor A. 14).

Since $\phi_n(D_n)$ converges to $D$ with respect to the Hausdorff distance, equation (\ref{oppositedrift}) guarantees $\IE \left( \tau^{\leftarrow}_{z,n,D_n} \right)$ is at least a power of $n$, for all $z \in \Zd$ such that $\left( \left(n^{-\frac{1}{d+1}}z_i\right)_{i\in \{1 \cdots d-1\}}, n^{-\frac{2}{d+1}} z_d \right) \in D_{\epsilon} \cap K.$ Hence, there is a constant $c_2$ uniform in $z$ such that:
$$ \IE \left( \mathcal{M}_{ n}(z) \right) \geq \IE \left( \tilde{\mathcal{L}}_{ n}(z) \right) \geq c_2 n^{\frac{2}{d+1}}.$$

Hence, both $ \IE \left( \tilde{\mathcal{L}}_{ n}(z) \right)$ and $\IE\left(\mathcal{M}_{ n}(z) \right) $ are at least powers of $n$.

As a consequence, with a probability that is exponentially close to $1$, uniformly in $z$ such that $\left( n^{-\frac{1}{d+1}}z_{1 \cdots d-1}, n^{-\frac{2}{d+1}} z_d \right) \in D_{\epsilon} \cap K,$ the point $z$ lies inside $A(n)$.

Using a simple union bound on the polynomial number of points $z$ such that: $\left( \left(n^{-\frac{1}{d+1}}z_i\right)_{i\in \{1 \cdots d-1\}}, n^{-\frac{2}{d+1}} z_d \right) \in D_{\epsilon} \cap K,$ the Borel-Cantelli lemma states that almost surely, for all but finitely many $n$,
$$D_{\epsilon} \subset \phi_n\left( A(n)  \right).$$



\begin{figure}
\includegraphics[scale=0.21]{fig2.png}
\includegraphics[scale=0.21]{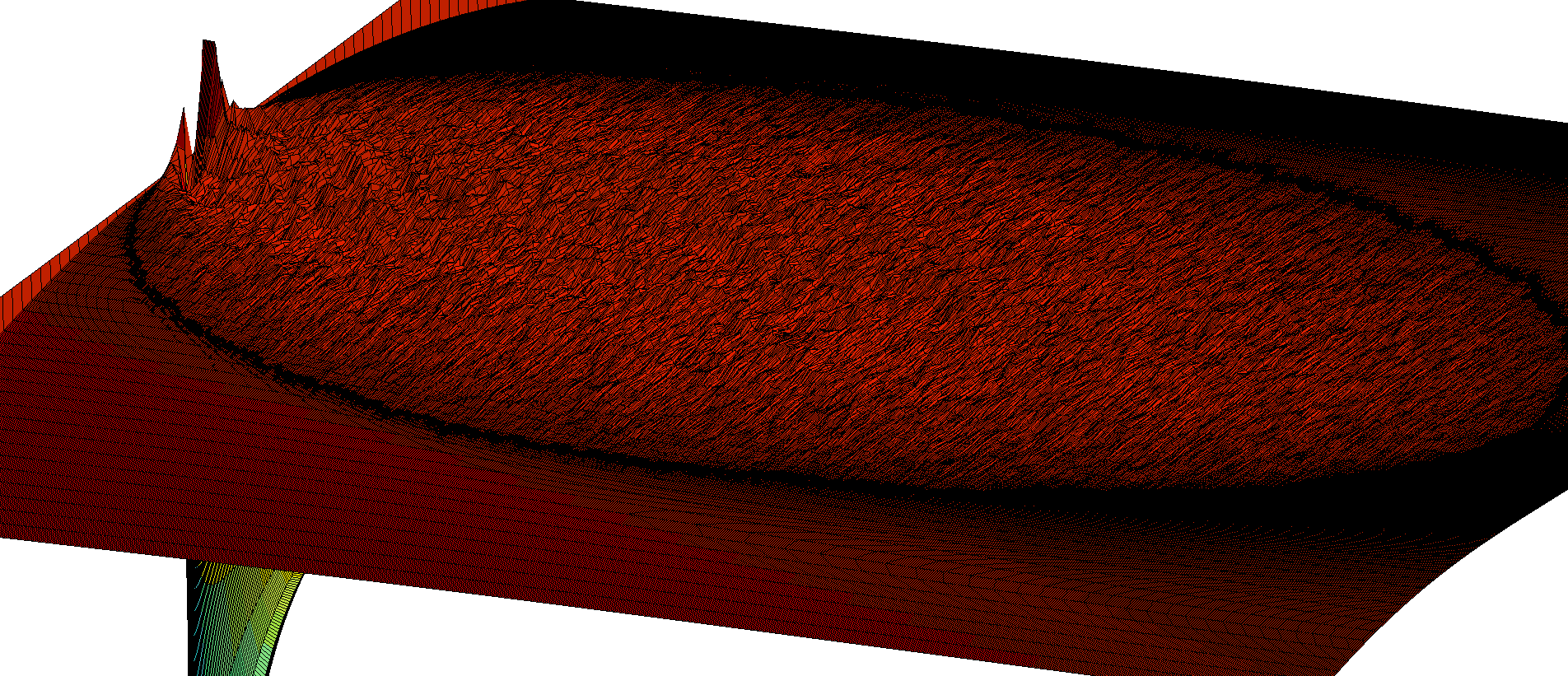}
\caption{On the top, the obstacle function $\gamma(x,t)$, and on the bottom, the sum of the obstacle function and normalised odometer. The drifted iDLA cluster can be glimpsed as the non-coincidence set.}\label{fig:3}
\end{figure}

\textbf{Remark:} The odometer function for the iDLA model does not appear directly in this proof. However, the function $f_{n,D_n}$ is closely linked to it. Hence, there is reason to believe that the normalised odometer for the drifted iDLA model also converges towards the difference between the obstacle function $\gamma$ and its least supercaloric majorant, as suggested by the simulation presented in Figure \ref{fig:3}.

\begin{theorem}
[External convergence towards $D$]
\label{extCV}
Let $A(n)$ be the drifted iDLA cluster.

Then, for every $\epsilon >0$, 
$$ \phi_n(A(n)) \subset D^{\epsilon},  $$
for all but finitely many $n$, almost surely.
\end{theorem}

First, we remark that the regularity of the solution to the obstacle problem is such that the Lebesgue measure of the boundary of $D$ is $0$, as stated in proposition \ref{regularity}. Together with Theorem \ref{idla-int}, it means that for every $\epsilon >0$, we can find $\eta >0$ and $K$ such that the number of points $z$ such that: 
$$\left( \left(n^{-\frac{1}{d+1}}z_i\right)_{i\in \{1 \cdots d-1\}}, n^{-\frac{2}{d+1}} z_d \right) \in D_{\epsilon} \cap K$$
is more than $(1 - \epsilon)n$.

Using the classical approach of \cite{lawler1992internal}, we argue that we only need control over the remaining $\epsilon n$ points. We will need to show that these particles do not make the cluster grow beyond distance $C\epsilon$ of its internal limiting shape $D$, with $C$ a constant depending only on the dimension. We will do this using a layers argument.

First, we need to argue that the particles spread out considerably when they travel a macroscopic distance in either time (the direction of the drift) or space (other dimensions). In the direction of the drift, this is ensured by the local central limit theorem. The case of other dimensions is studied in the following lemma.

For simplicity reasons, we begin by stating the desired result using a simple random walk in the $(d-1)$ space dimensions.

\begin{lemma}
[First contact point with a height set]

\label{ho}

Let $X_n$ be a $(d-1)$-dimensional simple random walk starting at the origin. Let $\tau_k$ be the hitting time of level $k$, namely $\tau_k = \inf\{ t\in \mathbb{N}, ||X_t||_{\infty} \geq k \}$.

Then, given $x_0>0$, there exists a constant $C$, such that for all $x,t \in \left( \mathbb{R}^*_+ \right)^2 $, with $x > x_0$, 
$$ \IP\left(\tau_{\lfloor x n^{\frac{1}{d+1}} \rfloor} = \lfloor t n^{\frac{2}{d+1}} \rfloor \right) \leq C n^{-\frac{d}{d+1}} .$$

\end{lemma}

This result can be explained as follows : the probability $P_n$ that $\lfloor t n^{\frac{2}{d+1}} \rfloor$ is the first hitting time of level $\lfloor x n^{\frac{1}{d+1}} \rfloor$ is the product of the probability that the value $\lfloor x n^{\frac{1}{d+1}} \rfloor$ is taken at time $\lfloor t n^{\frac{2}{d+1}} \rfloor$ by the probability, conditioned on this fact, that this level is never reached before:
\begin{eqnarray*}
P_n &=& \IP\left(\tau_{\lfloor x n^{\frac{1}{d+1}} \rfloor} = \lfloor t n^{\frac{2}{d+1}} \rfloor \right)\\
 &=& \IP \left( X_{\lfloor t n^{\frac{2}{d+1}} \rfloor} = \lfloor x n^{\frac{1}{d+1}} \rfloor \right) \IP \left( \forall t' < \lfloor t n^{\frac{2}{d+1}} \rfloor, ||X_{t'}||_{\infty} < \lfloor x n^{\frac{1}{d+1}} \rfloor \middle\vert X_{\lfloor t n^{\frac{2}{d+1}} \rfloor} = \lfloor x n^{\frac{1}{d+1}} \rfloor.  \right)
\end{eqnarray*}  

When $d=2$, a symmetry argument shows that the second factor is equal to the probability, for a walk starting at the origin and conditioned to reach $\lfloor x n^{\frac{1}{3}} \rfloor$ at time $\lfloor t n^{\frac{2}{3}} \rfloor$, to stay always strictly positive (except at the origin), and never go above $2\lfloor x n^{\frac{1}{d+1}} \rfloor$.

Considering only the first condition, our conditional probability is less than the probability, for a walk starting at the origin and conditioned to reach $\lfloor x n^{\frac{1}{3}} \rfloor$ at time $\lfloor t n^{\frac{2}{3}} \rfloor$, to stay always strictly positive (except at the origin).

This is a well-known problem (corresponding to the probability, when counting ballots,  that a candidate with $\lfloor xn^{\frac{1}{3}} \rfloor$ more votes than the other is always first on the tally when the total number of votes is $\lfloor tn^{\frac{2}{3}}\rfloor$), and this probability is exactly the following :
$$ \IP\left( \forall t' < t n^{\frac{2}{3}}, t'>0, X_{t'}>0 \middle\vert X_{t n^{\frac{2}{3}}} = x n^{\frac{1}{3}} \right) = \frac{\lfloor x n^{\frac{1}{3}} \rfloor}{\lfloor t n^{\frac{2}{3}}\rfloor} \leq C \frac{x}{t} n^{-\frac{1}{3}},$$
for a suitable constant $C >0$.

When $d >2$, we can apply the same argument. Indeed, suppose that the random walk reaches level $\lfloor x n^{\frac{1}{d+1}} \rfloor$ through coordinate $i$. Then consider the projection $X^i$ of $X$ to this coordinate. It is a lazy random walk, that remains in place with probability $\frac{d-2}{d-1}$ at each step and does a simple random walk otherwise.

For any $k$, the probability that $||X||_{\infty}$ stays less that $k$ can be bounded above by the probability that $|X^i|$ stays less than $k$.

Since $X^i$ is a lazy random walk, the number of non-zero steps it takes in time $\lfloor t n^{\frac{2}{d+1}} \rfloor$ is asymptotically greater than $A t n^{\frac{2}{d+1}}$, where $A$ is a suitable constant.

Hence, applying the previous result to $X^i$ yields :
\begin{eqnarray}
\IP \left( \forall t' < \lfloor t n^{\frac{2}{d+1}} \rfloor , ||X_{t'}||_{\infty} < \lfloor x n^{\frac{1}{d+1}} \rfloor \middle\vert  X_{\lfloor t n^{\frac{2}{d+1}} \rfloor} = \lfloor x n^{\frac{1}{d+1}} \rfloor \right) & \leq & \frac{\lfloor x n^{\frac{1}{d+1}} \rfloor}{A t n^{\frac{2}{d+1}}}, \nonumber\
\end{eqnarray}

Using the local central limit theorem, we get, for a suitable value of $A$ depending on $d$, (note that for $d=2$, one can take $A=1$): 
\begin{eqnarray}
\IP\left(\tau_{\lfloor x n^{\frac{1}{d+1}} \rfloor} = \lfloor t n^{\frac{2}{d+1}} \rfloor \right) &\leq&   \IP \left( X_{\rfloor t n^{\frac{2}{d+1}} \lfloor} = \rfloor x n^{\frac{1}{d+1}} \lfloor \right)  \frac{Cx}{At} n^{-\frac{1}{d+1}}\nonumber\\
 & \leq & \frac{C'x}{At} \exp \left( -\frac{||x||^2}{2t} \right) n^{-\frac{d-1}{d+1} -1} \nonumber\\
 & \leq & C'' n^{-\frac{d}{d+1}}, \nonumber
\end{eqnarray}
where the constant $C''$ is uniform in $x>x_0$. This concludes the proof of lemma \ref{ho}.

Let us now consider our drifted random walk. We consider its projection on the space dimensions, with a change of the time parameter such that is a simple random walk. Lemma \ref{ho} then states that the distribution of its crossing time for level $\lfloor xn^{\frac{1}{d+1}} \rfloor$ is spread out:
$$ \IP\left(\tau_{\lfloor x n^{\frac{1}{d+1}} \rfloor} = \lfloor t n^{\frac{2}{d+1}} \rfloor \right) \leq C n^{-\frac{d}{d+1}}.$$
We now argue our drifted random walk shares the same property, in the sense that the point at which its trajectory crosses the level $xn^{\frac{1}{d+1}}$ is spread out. Let us define $\tau_k = \inf\{t \in \mathbb{N}, \exists i \in \{1 \cdots d-1\} | \left(S(t)\right)_i | >k  \}$. Then there is a constant $C'$ uniform in $x>x_0>0$, $t$, such that:
$$ \IP\left( \left(S\left(\tau_{\lfloor x n^{\frac{1}{d+1}} \rfloor}\right) \right)_{d} = \lfloor t n^{\frac{2}{d+1}} \rfloor \right) \leq C' n^{-\frac{d}{d+1}}. \label{spreadho}$$
Indeed, the position in time of the crossing point is the number of steps taken by $S$ during the time it takes its projection to reach the level $xn^{\frac{1}{d+1}}$, which spreads out its distribution even more.




We will now bound the maximum distance that one of the $\epsilon n$ particles can reach outside the cluster in the space dimensions. 

We build the cluster and stop the particles when they exit the discrete version of $D_{\epsilon} \cap K$, so that we have a number $\alpha(n) \leq \epsilon n$ of particles waiting to be released on the boundary, and we define $\tilde A(i)$ as the cluster when $i$ of these particles have been released. It is a well-known property of the iDLA models that the law of $\tilde A(\alpha(n))$ does not depend on the order in which they are released. We will re-index these walks as $(S^{i})_{i \in \{1 \cdots \alpha(n)\}}$.


Choose $\eta_0 >0$, and, for $k \in \mathbb{N}$, consider the layer of points $$\mathcal{H}_k = \left\{ z\in \Zd, \lfloor n^{\frac{1}{d+1}}\eta_0 + k \rfloor \leq  d_h(z,D) \leq \lfloor n^{\frac{1}{d+1}} \eta_0 \rfloor + k +1 \right\},$$ where $d_h(z,D)$ denotes the distance, in $\mathbb{Z}^{d-1}$, between $z$ and set of points that normalise into $D$ in the $d-1$ non-drift directions (We ignore the distance in the direction of the drift). We will study the expected value of the number of particles in each layer $\mathcal{H}_k$. Le ut us define : 
$$\mu_{\mathcal{H},k}(l) = \IE\left( | \mathcal{H}_k \cap \tilde A(l) | \right).$$

When the $l+1$-th random walk $S^{l+1}$ is launched from the point at which it was stuck, consider the event that it adds to the cluster at a point of $\mathcal{H}_k$.
Now if $S^{l+1}$ adds to the cluster at a point of $\mathcal{H}_k$, it means that $S^{l+1}$ crosses $\mathcal{H}_{k-1}$ for the first time at a point of $\mathcal{H}_{k-1} \cap \tilde A(l)$. The distance that the random walk has to cross is bigger than $\eta_0 n^{\frac{1}{d+1}}$, so that using equation (\ref{spreadho}) gives a uniform bound the hitting probability of any given point. A simple union bound yields the following recursive relation, for a suitable constant $C_1$ :
\begin{eqnarray}
\mu_{\mathcal{H},k}(l+1) - \mu_k(l) & \leq & C_1 \mu_{\mathcal{H},k-1}(l) n^{-\frac{d}{d+1}} \nonumber
\end{eqnarray}

Summing over values $l$, and remarking that $\mu_{k-1}(0) = 0$, yields the following equation:
\begin{eqnarray}
\mu_{\mathcal{H},k}(\alpha(n)) & \leq & C_1 n^{-\frac{d}{d+1}} \sum_{l=1}^{\alpha(n)-1} \mu_{\mathcal{H},k-1}(l) \nonumber
\end{eqnarray}

By a simple induction argument, using the inequality $ \sum_{l=1}^{
\alpha(n)-1} l^{k-1} \leq \frac{\alpha(n)^k}{k}$, we have the following inequality:
\begin{eqnarray}
\mu_{\mathcal{H},k}(\alpha(n)) & \leq & \left( C_1 n^{-\frac{d}{d+1}} \right)^{k} \frac{\alpha(n)^{k}}{k!} \nonumber\\
         & \leq & \left( \frac{\alpha(n)}{k} C_1 n^{-\frac{d}{d+1}}  \right)^{k} \nonumber
\end{eqnarray}

Choosing $k = K \epsilon n^{\frac{1}{d+1}}$, with $K$ suitably big, results in the fact that $\mu_{\mathcal{H},k}(\alpha(n))$ is exponentially small. Hence, by the Borel-Cantelli lemma, the normalised cluster contains no points further away from $D$, in the direction of the drift, than $\eta_0+ K \epsilon$.

We proceed in the same way to bound the cluster in the direction of the drift.

Choose $\eta_0 >0$, and, for $k \in \mathbb{N}$, consider the layer of points $$\mathcal{V}_k = \left\{ z\in \Zd , \lfloor \eta_0 n^{\frac{1}{d+1}} \rfloor + k \leq d_v(z,D) \leq \lfloor \eta_0 n^{\frac{1}{d+1}} \rfloor +k+1 \right\},$$ where $d_v(z,D)$ denotes the distance, in $\mathbb{Z},$ between $z$ and points that normalise to $D$, in the direction of the drift.

We define $\mu_{\mathcal{V},k}(l)$ as the expected value of the number of particles settled in $\mathcal{V}_k$ after launching $l$ particles. 

Let us consider a walk $S^i$. Since all points of $\mathcal{V}_k$ are at least at distance $\eta_0n^{\frac{2}{d+1}}$ from the starting point of $S^i$, the local central limit theorem gives a uniform bound for the probability $\mathbb{P}_n$ of hitting any point of $\mathcal{V}_k$ before all the others: for a suitable constant $C_2$, we have:
$$\mathbb{P}_n \leq C_2 n^{\frac{d-1}{d+1}}.  \label{LTCL2}$$

Once more, since a particle needs to cross $\mathcal{V}_{k-1}$ on a point of $\tilde A(l)$ in order to contribute to $\tilde A(l+1)$ at a point of $\mathcal{V}_k$, we have a recursive relation that is similar to that of the previous proof, and for a suited constant $C_3$,
\begin{eqnarray}
\nu_{\mathcal{V},k}(l+1) - \nu_{\mathcal{V},k}(l) & \leq & C_3 \nu_{\mathcal{V},k-1}(l) n^{-\frac{d-1}{d+1}} \nonumber\\
\nu_{\mathcal{V},k}(\alpha(n)) & \leq & C_3 n^{-\frac{d-1}{d+1}} \sum_{l=1}^{\alpha(n)-1} \nu_{\mathcal{V},k-1}(l). \nonumber
\end{eqnarray}

By induction, this yields:
\begin{eqnarray}
\nu_{\mathcal{V},k}(\alpha(n)) & \leq & \left( C_3 n^{-\frac{d-1}{d+1}} \right)^{k} \frac{n^{k}}{k!} \nonumber\\
         & \leq & \left( \frac{\alpha(n)}{k} C_3 n^{-\frac{d-1}{d+1}}  \right)^{k}. \nonumber
\end{eqnarray}

Choosing $k = K \epsilon n^{\frac{2}{d+1}}$, with $K$ suitably big, yields that $\nu_{\mathcal{V},k}(\alpha(n))$ is exponentially small. Hence, by the Borel-Cantelli lemma, the normalised cluster contains no points at distance greater than $\eta_0 + K\epsilon $ from $D$. This concludes the proof of Theorem \ref{extCV}

\subsection{A more natural class of drifted random walks}

\label{extension}

Throughout our proofs so far, we have considered the drifted random walks $S^j$ with the following law:
\begin{eqnarray}
\IP\left( S^j (t+1) - S^j(t) = \pm e_i \right) &=& \frac{1-p}{2(d-1)} \qquad \text{ for $i=1 \cdots d-1,$ and} \nonumber\\
\IP\left( S^j (t+1) - S^j (t) = e_d \right) &=& p. \nonumber\\ 
\end{eqnarray}

A more natural drifted random walk can be defined as a walk that, at each time step, performs a simple random walk on $\Zd$ with probability $1-p$, and a step in the direction of the drift with probability $p$. Such a random walk $S'$ has the following law:
\begin{eqnarray}
\IP\left( S' (t+1) - S'(t) = \pm e_i \right) &=& \frac{1-p}{2d} \qquad \text{ for $i=1 \cdots d-1,$} \nonumber\\
\IP\left( S' (t+1) - S' (t) = \text{  }\text{  }e_d \right) &=& p+ \frac{1-p}{2d}, \text{ and} \nonumber\\
\IP\left( S' (t+1) - S' (t) = -e_d \right) &=& \frac{1-p}{2d}. \nonumber
\end{eqnarray}

In order to extend our limiting shape result to the cluster built using this random walk, we will see step by step how our proofs need to be modified.

First, we need a different version of the unfair divisible sandpile, in which the mass emitted from a point during its toppling is distributed in accordance with the new law of our random walk: a fraction $\frac{1-p}{2d}$ is sent to each neighbor, and an additional fraction $p$ is sent towards the neighbor in the direction of the drift.

Once again, we normalise our cluster by $n^{\frac{1}{d+1}}$ in non-drift directions, and $n^{\frac{2}{d+1}}$ in the direction of the drift.

Remark that the new version of the discrete operator we define still satisfies the minimum principle, so that the final configuration of mass is once again defined as the lowest supercaloric majorant of a suitable obstacle function, where the term supercaloric is defined with respect to the new discrete operator.

One suitable obstacle function is the following:
$$ \gamma_n(z) = -\frac{d}{d-1}\sum_{i=1}^{d-1} z_i^2 + z_d - n g(0,z), $$
where $g$ is the discrete Green function for the random walk $S'$.

Coupling our new drifted random walk with a random walk on $\mathbb{Z}^{d-1}$ with a suitable law, we get an convergence of the normalised obstacle function towards the function:
$$ \gamma'(x,t) = t - \frac{d}{d-1}|x|^2 - \frac{1}{p} \left( \frac{\beta}{\pi t} \right)^{\frac{d-1}{2}} \exp \left( - \beta \frac{||x||^2}{t} \right), $$
where $\beta = \frac{d p}{2(1-p)},$ that is once more uniform on compacts of $\IR^{d-1} \times  \IR^*_+$. We also have a suitable bound on the error.

The next step of the proof is to adapt the continuous operator $\mathfrak{K}$, which becomes:
$$ \mathfrak{K}(f) = \frac{1-p}{2d} \Delta f - p \frac{\partial f}{\partial t}. $$
Then, following the proof of lemma \ref{CVUu}, we can prove that the normalised version of the least supercaloric majorant of $\gamma_n$ converges to the least supercaloric majorant of $\gamma'$, where the supercaloric functions are defined as in section \ref{heatequation}, with respect to the new heat operator:
$$ \mathfrak{K}'f = \frac{1-p}{2d} \Delta f - p \frac{\partial f}{\partial t}.$$

This convergence is once again a little technical to prove, but it relies only on the precise estimates of the convergence of the Green function, and the fact that our discrete operator is a good approximation of the continuous one for functions with sufficient regularity. Remark that the terms corresponding to a simple random walk in the direction of the drift become negligible because of our normalisation.

Using the same arguments as in the proof of Theorem \ref{convergenceDS}, we can thus prove that the normalised cluster for the new unfair divisible sandpile model, when intersected with a compact, converges towards the limiting shape $D\cap K$, where $D$ is defined as the non-coincidence set of $\gamma'$ and its least supercaloric majorant.

The proof of the extension of the convergence result to the iDLA cluster follows that of Theorem \ref{idla-int} almost \textit{verbatim}, and the majoration of the exterior error can be derived from the proof of Theorem \ref{extCV}, using similar estimates on the probability for a random walk of hitting any given point in a set, after covering a macroscopic distance in any direction.







\section{Properties of the cluster}

We know that our cluster converges to a limiting shape $D$ defined in terms of an obstacle function, however we are interested in the properties of this limiting shape. To natural questions arise, the first being its boundedness. Indeed, our normalisation only truly captures the behaviour of the cluster if we can prove that the limiting shape under this normalisation is bounded. Moreover, we know that the limiting shape is a true heat ball, and the existence of a bounded true heat ball is a question of interest in PDE theory.

The second natural question concerning the shape $D$ is that of its universality: when we change the parameters of the model, we will show that the limiting shape is only modified through a simple variable change.

\subsection{Bounds on the cluster}
\label{bounds}

In this section, we work once again with our drifted random walk $S$, that performs at each step a simple random in $\Zd$ with probability $1-p$, and steps in the direction of the drift with probability $p$. Using probabilistic arguments yields bounds on the iDLA cluster, which can in turn be transferred to the continuous shape $D$. Hence, our limiting shape $D$ will be a bounded true heat ball.

\subsubsection{Horizontal bound}

\begin{lemma}
[Non-drift direction bound on the cluster]

The iDLA cluster $A(n)$, normalised by $n^{\frac{1}{d+1}}$, is bounded in all "non-drift" directions.
\label{lemmaH}
\end{lemma}
Proof:
We proceed as in the proof of Theorem \ref{extCV}, and we consider the intersection of the entire cluster with strips in successively non-drift and drift directions.

Choose $x_0 >0$, and define $k_0 = \lfloor x_0 n^{\frac{1}{d+1}} \rfloor$, and $\Gamma_k = \{ (x,t) \in \mathbb{R}^{d-1} \times \mathbb{R},  ||x||_{\infty} = k_0+k \}$
In this proof, we study the expected value of the number of particles in each cylindrical domain $\Gamma_k$.
Define $\mu_k(l) = \IE\left( \left| A(l) \cap \Gamma_{k}  \right| \right) $.

When the $l+1$-th random walk $X^{l+1}$ is launched from the origin, consider the event that it adds to the cluster at a point of $\Gamma_k$.
Now if $X^{l+1}$ adds to the cluster at a point of $\Gamma_k$, it means that $X^{l+1}$ crosses $\Gamma_{k-1}$ for the first time at a point of $\Gamma_k \cap A(l)$.
Lemma \ref{spreadho} gives a uniform bound on this hitting probability, so that we have the following recursive relation, for a suitable constant $C_1$ :
\begin{eqnarray}
\mu_k(l+1) - \mu_k(l) & \leq & C_1 \mu_{k-1}(l) n^{-\frac{d}{d+1}} \nonumber
\end{eqnarray}

Summing over values $l$, and remarking that $\mu_{k-1}(0) = 0$, yields the following equation:
\begin{eqnarray}
\mu_k(n) & \leq & C_1 n^{-\frac{d}{d+1}} \sum_{l=1}^{n-1} \mu_{k-1}(l) \nonumber
\end{eqnarray}

By a simple induction argument, using the inequality $ \sum_{l=1}^{n-1} l^{k-1} \leq \frac{n^k}{k}$, we have the following inequality:
\begin{eqnarray}
\mu_k(n) & \leq & \left( C_1 n^{-\frac{d}{d+1}} \right)^{k} \frac{n^{k}}{k!} \nonumber\\
         & \leq & \left( \frac{n}{k} C_1 n^{-\frac{d}{d+1}}  \right)^{k} \nonumber
\end{eqnarray}

Choosing $k = x_1 n^{\frac{1}{d+1}}$, with $x_1$ suitably big, results in the fact that $\mu_k(n)$ is exponentially small. Hence, by the Borel-Cantelli lemma, the normalised cluster is almost surely asymptotically bounded by $x_0+x_1$.

\subsubsection{Vertical bound}
\begin{lemma}
[Drift-direction bound on the cluster]
The iDLA cluster normalised by $n^\frac{2}{d+1}$ in the direction of the drift is bounded in the direction of the drift.
\label{lemmaV}
\end{lemma}

Proof: The Local Central Limit theorem gives a uniform bound for the hitting probability of a vertical strip at a particular point: for a suitable constant $C_2$, we have:
$$\IP\left( X_{\lfloor tn^{\frac{2}{d+1}} \rfloor } = \lfloor x n^{\frac{1}{d+1}} \rfloor \right) \leq C_2 n^{\frac{d-1}{d+1}}  \label{LTCL}$$

We will follow a similar strategy as in the previous proof.

Choose $t_0 >0$, and define $k_0 = \rfloor t_0 n^{\frac{2}{d+1}} \lfloor $ and the vertical strips $V_k = \{ (x,t) \{ (x,t) \in \mathbb{R}^{d-1} \times \mathbb{R}, t = k_0+  k \}$.

We will study the number of particles in each vertical strip $V_k$. Let $\nu_k(l) = \IE\left( \left| A(l) \cap C_{k_0+k} \right| \right) $ denote its expected value.

Since a particle needs to cross $V_{k-1}$ on a point of $A(l)$ in order to contribute to $A(l+1)$ at a point of $V_k$, we have a recursive relation that is similar to that of the previous proof, and for a suited constant $C_3$,
\begin{eqnarray}
\nu_k(l+1) - \nu_k(l) & \leq & C_3 \nu_{k-1}(l) n^{-\frac{d-1}{d+1}} \nonumber\\
\nu_k(n) & \leq & C_3 n^{-\frac{d-1}{d+1}} \sum_{l=1}^{n-1} \nu_{k-1}(l) \nonumber
\end{eqnarray}

By induction, this yields:
\begin{eqnarray}
\nu_k(n) & \leq & \left( C_3 n^{-\frac{d-1}{d+1}} \right)^{k} \frac{n^{k}}{k!} \nonumber\\
         & \leq & \left( \frac{n}{k} C_3 n^{-\frac{d-1}{d+1}}  \right)^{k} \nonumber
\end{eqnarray}

Choosing $k = t_1 n^{\frac{2}{d+1}}$, with $t_1$ suitably big, yields the summability of $\nu_k(n)$. Hence, by the Borel-Cantelli lemma, the normalised cluster is almost surely asymptotically bounded by $t_0+t_1$.

\subsection{Rescaling}

In this section, we will show that the various limiting shapes can be written in terms of one another by rescaling differently in the direction of the drift and in the other directions.
\begin{lemma}

Let $p_1$ and $p_2$ be two drift parameters in $\left( 0, \infty \right)$, and let $D_1$ (respectively $D_2$) be the limiting shape of the unfair divisible sandpile model run with parameter $p_1$ (respectively $p_2$). 

Then $D_2$ is the image of $D_1$ by a change of variables $x \rightarrow \mu x, t \rightarrow \lambda t$.
\end{lemma}

First, we introduce the additional parameter $k\in \IR^*_+$, which measures the quantity of mass sent from the origin. We will call unfair divisible model with mass $k$ the unfair divisible sandpile model run with initial mass $k n$ at the origin. Remark that the limiting shape of this model is of course obtained from that of the original model by a rescaling (with different coefficients in the drift and non-drift directions). Moreover, the limiting shape of the model can still be obtained as the solution to our parabolic obstacle problem, the only difference being an adjusted coefficient $k$ in the last term of the obstacle function:

\begin{eqnarray}
\gamma_k(x,t) &=& t-||x||^2 - \frac{k}{p}  \left( \frac{\beta}{\pi t} \right)^{\frac{d-1}{2}} \exp \left( -\beta \frac{||x||^2}{t} \right). \nonumber
\end{eqnarray}

We are now ready to compare $D_1$ and $D_2$. Since the shape $D$ stems from the obstacle $\gamma$, we only need to show that $\gamma_2$, the obstacle function for parameter $p_2$ can be related to $\gamma_{1,k}$, the obstacle function for drift $p_1$ and mass $k$. Consider the two following functions:
\begin{eqnarray}
\frac{1}{C} \gamma_{1,k} (\mu x, \lambda t) & = & \frac{\lambda}{C} t - \frac{\mu^2}{C}||x||^2 - \frac{k}{C p_1  \lambda^{\frac{d-1}{2}}}\left( \frac{\beta_1}{\pi t} \right)^{\frac{d-1}{2}} \exp \left( -\beta_1 \frac{\lambda}{\mu^2} \frac{||x||^2}{t} \right), \nonumber\\
\gamma_2(x,t) &=& t-||x||^2 - \frac{1}{p_2}  \left( \frac{\beta_2}{\pi t} \right)^{\frac{d-1}{2}} \exp \left( -\beta_2 \frac{||x||^2}{t} \right), \nonumber
\end{eqnarray}
with $C$ a suitable constant such that $$\mathfrak{K} \left(\frac{\lambda}{C} t - \frac{\mu^2}{C} |||x||^2 \right) = -1.$$


We see that the constant $\beta_1$ can be changed to $ \beta_2=\frac{\lambda}{\mu^2}\beta_1$, provided that $k$ is such that the last term of $\frac{1}{C}\gamma_{1,k}$ is equal to that of $\gamma_2$, that is to say, 
$$\frac{k_1}{C p_1  \lambda^{\frac{d-1}{2}}} = \frac{1}{p_2}.$$
Since  $\frac{1}{C} \gamma_{1,k} (\mu x, \lambda t)$ and $\gamma_2(x,t)$ now only differ by a caloric function, they give rise to the same limiting odometer, hence to the same limiting shape. We conclude using the first part of the proof to recover that $D_1$ and $D_2$ are indeed images of one another  by a transformation of the required form $x \rightarrow \mu x, t \rightarrow \lambda t$.

\section{Conclusion}

At this point, we have proved the convergence of our two models towards a limiting shape $S$ that solves the following PDE problem:
given $\phi$ a $C^{\infty}$ function of time and space such that
$$ \frac{1-p}{2(d-1)} \Delta \phi + p \frac{\partial \phi}{ \partial t} =0. $$
Then the following mean value property holds:
$$ \int_{S} \phi(z,t) dzdt = |S| \phi(0). \label{trueheatball}$$

 We have also proved, using probabilistic estimates on random walks, that the iDLA cluster is bounded. Moreover, the regularity of the problem enables us to transpose these bounds to $S$, since $S$ is sufficiently defined in terms of a Hausdorff limit. This concludes the proof of Theorem \ref{maxitheoreme}.

\end{document}